\documentclass[12pt]{amsart}
 
\usepackage{latexsym, amsfonts, amsmath, amssymb, amscd, epsfig, xypic} 
 
\newtheorem{theorem}{Theorem}[section]
 
\newtheorem{corollary}[theorem]{Corollary}
 
\newtheorem{lemma}[theorem]{Lemma}

\theoremstyle{definition}
 
\newtheorem{definition}[theorem]{Definition}

\newtheorem{remark}[theorem]{Remark}

\newtheorem{example}[theorem]{Example}
 
\theoremstyle{parrafo}
 
\newtheorem{parrafo}[theorem]{{\!}}
 
\numberwithin{equation}{theorem}
 
\newcommand{\call}{\mathcal{L}}

\newcommand{\nat}{\mathbb N}
 
\newcommand{\rac}{\mathbb Q}
 
\newcommand{\ent}{\mathbb Z}
 
\newcommand{\ma}{\mathbb A}

\newcommand{\cali}{{\mathcal {I}}}

\newcommand{\cala}{{\mathcal {A}}}
 
\newcommand{\calo}{{\mathcal {O}}}

\newcommand{\Reg}{\mbox{Reg }}
 
\newcommand{\Sing}{\mbox{Sing}}
 
\newcommand{\RSing}{\mbox{RSing}}

\newcommand{\Max}{\mbox{\underline{Max} }}
 
\newcommand{\word}{{w-ord}}

\newcommand{\wordid}{{w-ord}^d_{i}\mbox{ }}

\newcommand{\tqz}{({\mathbb Q}\times {\mathbb Z},\leq )}
 
\begin{document}

\title[]{A Strong Desingularization Theorem}
 
\author{A. Bravo and O. Villamayor U.}
 
\address{Dpto. Matem\'aticas, Facultad de Ciencias, Universidad
Aut\'onoma de Madrid, Canto Blanco 28049 Madrid, Spain}
\email{ana.bravo@uam.es; villamayor@uam.es}
 
\address{Current address of Ana Bravo: Department of Mathematics, University of Michigan, Ann Arbor, MI 48109-1109}
\email{anabz@math.lsa.umich.edu} 
 
\thanks{Supported by DGICYT 1317 DP Spain}
 
 
 
\date{}
 
\dedicatory{}
 
\commby{}

\begin{abstract}
Let $X$ be a closed  subscheme embedded in a
 scheme $W$, smooth over a field ${\bf k}$ of characteristic zero,
and let $\cali(X)$ be the sheaf of ideals defining $X$. Assume that 
the set of regular points of $X$ is dense in $X$. 
We prove that there exists a proper, birational morphism,
$\pi: W_r\longrightarrow W$, obtained
 as a composition of monoidal transformations,
  so that if $X_r\subset W_r$ denotes the strict transform of $X\subset W$
then:

1) The morphism $\pi:W_r\longrightarrow W$ is an embedded
 desingularization of $X$ (as in Hironaka's Theorem);

2) The {\em total transform} of $\cali(X)$ in $\calo_{W_r}$ factors as
a product
of an invertible sheaf of ideals $\call$ supported on the
exceptional locus,
and  the sheaf of ideals defining the strict transform of $X$ (i.e.
$\cali(X)\calo_{W_r}=\call\cdot\cali(X_r)$).

This result is stronger than Hironaka's Theorem, in fact (2) is novel and does not hold for  desingularizations which follow Hironaka's line of proof unless $X$ is a hypersurface.
 We 
will say that $W_r\longrightarrow W$ defines a {\em Strong
Desingularization of $X$}. 
\end{abstract}

\maketitle
 

\begin{center}
{\sc Introduction}
\end{center}
\markboth{A. Bravo and O. Villamayor U.}{A Strong Desingularization Theorem}
\vspace{0.5cm}

After  Hironaka's monumental work,  where a non-constructive,
existential proof of resolution of singularities
was given  (cf. \cite{Hir64}), several alternative 
proofs have appeared. Some of them are simplified weak non-constructive versions 
(cf. \cite{AbrJong97}, \cite{AbrWang97}, \cite{BogPan96}), 
 while others are  constructive algorithms of resolution of singularities   
(cf. \cite{BM97}, \cite{Villa89}, \cite{Villa92},
\cite{EncVil98}). The latest ones  provide proofs
of
desingularization
within Hironaka's line of development.

\

A simple alternative proof, conceptually different from Hironaka's,
was pointed out in an addendum in \cite{EncVil97:Tirol}
(see also \cite{EncVil99}).
There desingularization is achieved as a
byproduct of a much simpler problem, namely that of {\em principalization of ideals}.
Hilbert Samuel functions and normal flatness are avoided
in the proof of desingularization with this  strategy.

\

The alternative proofs of  \cite{EncVil97:Tirol} and \cite{EncVil99}
rely on the fact
that both  desingularization and principalization are special
 cases of the so called
 {\em resolution of basic objects} (cf. Section \ref{secbo}).
 Recently G.~Bodn\'{a}r and J.~Schicho
have implemented an algorithm of resolution of
basic
objects  in MAPLE,  providing thereby a computer program  of
embedded desingularization (cf. \cite{GabSch98} and \cite{GabSch99}). 

\

Within the framework of this new generation of algorithms 
where principalization and desingularization 
follow from suitable defined resolution of basic objects, 
 we obtain the new result presented here. By using a
different {\em algorithm of desingularization}, the form 
of desingularization we achieve is
 stronger than Hironaka's and the ones given by
 previous algorithms. Our result  is stated as follows:

\

\noindent {\bf Theorem of Strong Embedded 
Desingularization.} {\em Let $W$ be a smooth scheme over a field of characteristic 
zero ${k}$, and let $X\subset W$ be a closed subscheme. Assume that $\mbox{Reg}(X)$, the 
set of regular points in $X$, is dense in $X$  (e.g. if $X$ is reduced). Then there exists a proper 
birational morphism obtained as a composition of monoidal transformations,
$$\Pi_r:W_{r}\longrightarrow W,$$
so that if \( X_{r}\subset W_{r} \) is 
    the strict transform of \(X\), and if $E_r$ is is the exceptional locus of $\Pi_r$, then 
 locally at every point $y\in W_r$   there will be a factorization of the total transform 
of the form:
\begin{equation} 
\label{ceros}
J_{\pi(y)}\calo_{W_r,y}=\call_y(J^\prime)_y\subset \calo_{W_r,y}
\end{equation}
where $\call_y$ is principal and supported on the exceptional locus of $\pi$,
$J^\prime_y$ is  the weak transform of $J$ at $y$, and  the following conditions hold:  
\begin{enumerate}
\item[a)] If $J_{\pi(y)}=\calo_{W,\pi(y)}$, or if $J_{\pi(y)}\neq \calo_{W,{\pi(y)}}$ 
and $J_{\pi(y)}$  is a regular prime ideal, then $\calo_{W,{\pi(y)}}=
 \calo_{W_r,y}$.
\item[b)] If $J_{\pi(y)}\neq \calo_{W,_{\pi(y)}}$ and $J_{\pi(y)}$  is 
not a regular prime ideal, then there is a regular system of
parameters $\{x_1,\ldots,x_d\}\subset \calo_{W_r,y}$ and an integer $e\geq 0$ such that:
\begin{equation}
\label{hada1}
\call_y=x_{e+1}^{c_{e+1}}\cdot\ldots\cdot x_k^{c_k}, \mbox{ } \mbox{ with }\mbox{ } e<k\leq d
\end{equation}
and 
\begin{equation}
\label{hada}
J^\prime_y=<x_1,\ldots,x_e>, 
\end{equation}
and $J^\prime_y$ defines the strict transform of $X$ at $y$.
\end{enumerate}
}

\

We
will say that $W_r\longrightarrow W$ defines a {\em Strong
Desingularization of $X$}.

\

An announcement of this result appears in \cite{BV}. 

\

In Section \ref{lang} we will  reformulate our theorem  
to clarify how it relates to the classical formulation of desingularization  (see 
Theorem \ref{MainTh}). 
The fact that older results do not yield factorizations with the local conditions  
stated in (b)   will be  seen by  analyzing  the
 example where $W=\ma_{\rac}^3=\mbox{Spec }{\mathbb
{Q}}[x_1,x_2,x_3]$ and $X$ is the irreducible curve defined by
the ideal $\cali(X)=<x_1,x_2x_3+x_2^3+x_3^3>$ (see Section \ref{totaltransform} and 
 Example \ref{ejemplo1}). 
Our Strong Desingularization Theorem
 is therefore a  novel
algebraic formulation of embedded desingularization, which
extends the classical
result about  embedded hypersurfaces to
subschemes of arbitrary codimension (note that if $X$ is a hypersurface and if $\Pi_r:W_r\longrightarrow W$ 
defines an embedded resolution of $X$, then the local condition 
of part (b) always holds). 

\

In our development we do not consider the notion of 
strict transforms of ideals. Given an ideal defining a reduced subscheme, and a  
sequence of monoidal transformations, we study the primary 
decomposition of its total transform, and  we show that we can define a  
 sequence of monoidal transformations  so that all the new (exceptional) 
primary components are locally principal (defining the invertible ideal 
$\call$ which appears in (\ref{ceros})). 
Our result answers a questions formulated by A. Nobile, 
which was the starting point of this research. 

\

If $X\subset W$ is a complete intersection then a resolution of the
sheaf $\calo_X$ by free $\calo_W-$modules can be defined in terms
of a Koszul complex. A byproduct of our theorem is that such Koszul
complex (resolution) can be lifted to the Koszul complex providing
a {\em resolution} of $\calo_{X_r}$ (as $\calo_{W_r}-$module),
after a simple twisting by the invertible sheaf $\call$ appearing
in (iii) of Theorem \ref{MainTh} (cf. Corollary \ref{Koszul}).

\

{\em Idea of the proof of the Theorem of Strong Embedded Desingularization:} 
In this paper the approach to desingularization is substantially different from  that in earlier works. 
Let $J=\cali(X)$ be as in the theorem. Even if  $J$ is 
 not be reduced,  we are requiring  that  there is a 
{\em dense} open set of regular points in $X$. This is the traditional requirement on a subscheme 
in order to formulate the problem of desingularization. As a consequence 
 the primary components of $J$  
corresponding to the minimal associated prime ideals should be reduced. 
Now we turn to the question on how  the primary
decomposition of $J$ can be modified after a sequence of monoidal transformations on $W$. 
Given a sequence of monoidal transformations
$$
W_k\longrightarrow W_{r-1}\longrightarrow \ldots\longrightarrow W_1\longrightarrow W_0=W$$ 
we will consider the factorization, 
$$J{\mathcal O}_{W_i}=\call_i J^{\prime}_i, \mbox{ } \mbox{ for  }\mbox{ } i=1,\ldots,k$$ 
where $J^{\prime}_i$ denotes the weak transform of $J$ in $W_i$ 
and $\call_i$ is locally a monomial ideal  (see (\ref{ceros})).
The strategy of our proof is to show  first that we can define 
 a sequence of monoidal transformations 
$$
W_l\longrightarrow W_{l-1}\longrightarrow \ldots\longrightarrow W_0,$$
so that 
$V((J^{\prime}_l)_y)$ is locally included in a smooth hypersurface for every $y\in V(J^{\prime}_l)$. 
Then we  show
that by applying more monoidal transformations 
$$
W_m\longrightarrow W_{m-1}\longrightarrow \ldots\longrightarrow W_l,$$
we may require $V((J^{\prime}_m)_y)$
to be locally included in a smooth subscheme of codimension 2 for every $y\in V(J^{\prime}_m)$, and so on.
 This leads us to the notion of {\em embedded relative codimension} 
stated  in Definition \ref{omni}. 
This an algebraic local  condition on the factorization in (\ref{ceros}). If the
embedded relative codimension  is at least one, then $V((J^{\prime}_i)_y)$ is locally
included in a smooth subscheme of codimension 1. If the   embedded 
relative codimension
 is at least two, then $V((J^{\prime}_i)_y)$  is locally included in a smooth
subscheme of codimension two, and so on.
 Lemma \ref{lema}  shows that we can define a sequence of monoidal
transformations  so that the locally embedded codimension of the weak
transform coincides with the codimension of $X$ in $W$. 
It is in the proof of   this lemma where we make use of some properties of
algorithmic principalization. Following this argument 
we provide already a  simple geometric 
proof of Theorem \ref{MainTh} in \ref{proofoftheo}  .

\

The paper is organized as follows: After the formulation of Theorem \ref{MainTh} 
in Section \ref{lang}, Section \ref{totaltransform} is devoted  to explain and illustrate 
what makes this work different from previous approaches to desingularization. 
 In Section \ref{secbo} we  review the notion of {\em basic object} and that of  
{\em resolution of basic objects}, which is closely related 
to that of principalization of ideals (see Remark \ref{so}). As pointed out above, 
principalization of ideals plays an important role in the proof of Theorem \ref{MainTh}. 
In Section \ref{sketch} we show how   Theorem \ref{MainTh}
  follows from Lemma \ref{lema} and  a simple geometric argument (see \ref{proofoftheo}). 
 In Section \ref{last1}  we explain the relation between Lemma \ref{lema} 
 and  the two main invariants of 
algorithmic principalization (in particular Remark \ref{svive} 
is already a hint for the proof of Lemma \ref{lema}). This relation  
 is also illustrated with examples in Section \ref{examples}. Lemma \ref{lema} is finally 
proved in 
section \ref{last}. 

\

We are indebted to Prof. J.M. Aldaz for several useful suggestions during the 
preparation of this paper. 

\section{Formulation of the Main Theorem}
\label{lang} 
We briefly explain the notions of {\em pairs} and {\em
transformation of pairs} which are suitable for the
formulation of both {\em Strong Embedded Desingularization} and 
{\em Principalization of Ideals} (cf. Theorem \ref{MainTh} and
Definition \ref{princip} below).
 
\begin{definition}
\label{pairs}
\cite{EncVil97:Tirol}.
Let \( W \)  be  a pure
dimensional scheme, smooth over a field \( \mathbf{k} \) of
characteristic zero, and let \( E=\{H_{1},\ldots,H_{r}\} \) be a set of
smooth hypersurfaces in \( W \) with  normal crossings (i.~e.
\( \cup_{i=1}^{r} H_{i} \) has normal crossings). The couple  \( (W,E) \)
is said to be a {\em pair}.
A regular closed subscheme \( Y\subset W \) is said to be
{\em permissible} for
the pair \( (W,E) \) if \( Y \) has
normal crossings with  \( E \).

If \( Y \subset W\) is permissible  for a pair $(W,E)$, we
define a {\em transformation  of  pairs} in the following way:
Consider the blowing-up with center $Y$,
 $W\stackrel{\Pi}{\longleftarrow} W_1$,
and define \( E_{1}=\{H'_{1},\ldots,H'_{r},H_{r+1}\} \), where
\(H'_{i} \) denotes the strict transform of \( H_{i} \), and \(
H_{r+1}$ denotes $\Pi^{-1}(Y) \) the exceptional hypersurface in \( W_{1} \).
Note that \(W_{1}\) is smooth and that \(E_{1}\) has normal
crossings. We say that $(W,E)\longleftarrow(W_{1},E_{1})$ is a
{\em transformation of the pair $(W,E)$. }
\end{definition}

\begin{parrafo}
\label{MainTh} {\rm {\bf Main Theorem }}(of Strong Embedded
Desingularization) Let $(W_0, E_0=\{\emptyset\})$ be a pair and let
$X_0\subset W_0$
    be  a closed subscheme defined by
    $\cali(X_0) \subset \calo_{W_0}$. Assume that the open set \( \Reg(X) \)
of regular points is
    dense in $X$ (e.g. if $X$ is reduced).
 Then there exists a
    finite sequence of transformations of pairs
\begin{equation}
\label{aparece}
(W_0,E_0)\longleftarrow \cdots \longleftarrow (W_r,E_r)
\end{equation}
    inducing a proper birational morphism
     $\Pi_r:W_{r}\longrightarrow W_{0}$,
so that
    setting $E_r=\{H_1,\ldots,H_r\}$ and letting \( X_{r}\subset W_{r} \) be 
    the strict transform of \(X_0\), we have that: 
\begin{enumerate}
\item[(i)] $X_r$ is regular in $W_r$, and $W_r\setminus\cup_{i=1}^rH_i\simeq
W_0-\mbox{Sing}(X)$. In particular
$\Reg(X)\cong\Pi^{-1}_r(\Reg(X))\subset X_{r}$ via \( \Pi_r \).
\item[(ii)] \( X_{r} \) has normal crossings with \( E_{r}=\cup_{i=1}^rH_i
\) (the exceptional locus of $\Pi_r$).
\item[(iii)] The total transform of the ideal $\cali(X_0)\subset \calo_{W_0}$
 factors as a product of ideals in $\calo_{W_r}$:
$$\cali(X)\calo_{W_r}={\mathcal {L}}\cdot\cali(X_r),$$
where now $\cali(X_r)\subset \calo_{W_r}$ denotes the sheaf of ideals defining
$X_r$, and ${\mathcal {L}}=\cali(H_1)^{a_1}
\cdot\ldots\cdot\cali(H_r)^{a_r}$ is an invertible sheaf of ideals
supported on the exceptional locus of $\Pi_r$.
\end{enumerate}
\end{parrafo}
 
\

Note that part (i) of Theorem \ref{MainTh} ensures that the only points of $W_0$ that will be modified 
by the morphism $\Pi_r$
 are the ones in $\Sing(X_0)$. 
In fact, the proof of the Theorem \ref{MainTh}, 
   follows from the simpler problem of principalization.

\begin{definition}
\label{princip} Let $I\subset \calo_W$ be a sheaf of ideals. A {\em
principalization of $I$} is a proper birational morphism
$W_1\longrightarrow W$ such that   $I\calo_{W_1}$ is
an invertible
sheaf of ideals in $W_1$. A {\em strong principalization of $I$} is a chain
of transformations of pairs
$$(W_0,E_0=\emptyset)=(W,E)\longleftarrow \ldots\longleftarrow (W_r,E_r)$$
such that $W\longleftarrow W_r$ defines an isomorphism over the open subset
$W\setminus
V(I)$, and
$${\mathcal {L}}=I\calo_{W_r}=\cali(H_1)^{c_1}
\cdot\ldots\cdot\cali(H_s)^{c_s},$$
where $E^{\prime}=\{ H_1,H_2,\dots, H_s\}$ are regular hypersurfaces with normal
crossings and for $i=1,\ldots,s$,  $c_i\geq 1$. If  $V(I)$ is of
codimension $\geq 2$, this means that $E^{\prime}=E_r$ and the total
transform of $I$ is locally spanned by a monomial supported on the
exceptional locus of
$\Pi_r: W_r \to W$.
\end{definition}

\section{Total transform versus strict transform}
\label{totaltransform}
Our next goal is to explain the difference between the algorithm 
of desingularization that we present in Theorem \ref{MainTh} and 
the previous ones which follow from Hironaka's line of proof. 
To do so we will study the example of desingularization of a curve 
embedded in a three-dimensional space. This example together with the discussion that 
follows it already appear in \cite[Section 3]{BV}. We reproduce it 
here for the convenience of the reader.

\

 Let $W_0=A^3_{\rac}=\mbox{Spec}({\rac}[x_1,x_2, x_3])$  and consider the
curve $C$
defined  by
$$\cali(C) =<x_1,x_2x_3+x_2^3+x_3^3>.$$ 
Let  $W_0\longleftarrow W_1$ be  the quadratic transformation
at the origin, $H\subset W_1$ the exceptional divisor, and $C_1$ the
{\em strict transform} of $C$.
This defines an embedded desingularization of $C$, in the usual sense,
since both (i) and (ii) of Theorem \ref{MainTh} hold.
 
\
 
{\bf A)} {\em (On condition \ref{MainTh} (iii))}. Since the ideal $\cali(C)$ has order 1
at the
center of the quadratic transformation,
the {\em total transform}  of $\cali(C)$, namely $\cali(C){\mathcal
O}_{W_1}$, can be factored as a product,
$\cali(C){\mathcal O}_{W_1}=\cali(H)^1\overline{J}_1$ for some coherent ideal
$\overline{J}_1\subset {\mathcal O}_{W_1}$ which {\em does
not } vanish along $H$.
 
\ 
 
Note that $\cali(C_1)$ is a primary component of $\overline{J}_1$.
However, $\overline{J}_1\subsetneqq\cali(C_1)$, and hence Theorem \ref{MainTh} (iii) does not hold.
To see why, it is convenient to express both ideals in terms of conductors: By
definition,
$$\overline{J}_1=(\cali(C){\mathcal O}_{W_1}:\cali(H)^1).$$
On the other hand, the ideal of the strict transform is an increasing union 
$$\cali(C_1)=\cup_{k\geq 0}(\cali(C){\mathcal O}_{W_1}:\cali(H)^k),$$
or, in other words, $\cali(C_1)=(\cali(C){\mathcal O}_{W_1}:\cali(H)^N)$
for $N$ large enough, since
$$(\cali(C){\mathcal O}_{W_1}:\cali(H)^k)\subset (\cali(C){\mathcal
O}_{W_1}:\cali(H)^{k+1}).$$
 
In this example $H\simeq P^2_{\rac}$ and $C_1$ cuts ${\mathbb P}^2_{\rac}$ transversally at two
different
points. Let $L\subset {\mathbb P}^2_R$ be
the line defined by these two points, and let $\cali(L)\subset {\mathcal
O}_{W_1}$ be the ideal of $L (\subset W_1)$. Then 
$$\overline{J}_1=(\cali(C){\mathcal
O}_{W_1}:\cali(H))\subsetneqq\cali(C_1)=(\cali(C){\mathcal
O}_{W_1}:\cali(H)^2),$$
and looking at a suitable affine chart it follows that $\cali(L)$ is a
primary component of $\overline{J}_1$, but  (of course), not
of $\cali(C_1)$ (see Example \ref{ejemplo1} where the same example is 
 studied with more detail).
Therefore (iii) of
Theorem \ref{MainTh} does not hold for the embedded desingularization
defined by $\Pi$.
 
\
 
In Hironaka's line of proof the centers of monoidal
transformations, chosen in accordance to the so called {\em standard basis}, are always
{\em included in the strict transform} of the scheme, so only {\em quadratic transformations} are 
applied in the one dimensional case. In the case of our
singular curve, the first monoidal transformation must be the one we have
defined above, and any other center will have dimension zero. Now
 $\cali(L)$ is a primary component of $\overline{J}_1$ supported on $L$, which has
dimension 1; so we will never eliminate $\cali(L)$ as a primary component of the 
total transform of $\cali(C)$ by blowing up
closed points; hence (iii) will never hold for  desingularizations
of this curve that follow from  Hironaka's proof.
 
\
 
In order to achieve (iii) one must blow up $L$ (or some strict transform
of $L$). Using the new algorithm
that we propose, we first consider the quadratic transformation $\Pi:
W_1\longrightarrow W_0$, and
then we  blow-up at the one dimensional scheme $L$. Since $L\subset H$, the first isomorphism in
Theorem \ref{MainTh} (i) is preserved after such monoidal
transformation.
 
\
 
{\bf B)} {\em (On a question of complexity).} We think of a subscheme $X$
of a smooth
scheme $W$, at least locally, as a finite number of {\em equations}
defining the ideal $\cali(X)$. An  algorithm of
desingularization should  provide us with:
\begin{enumerate}
\item[(1)] A sequence of monoidal transformations over the smooth scheme $W$,
say $W_n\longrightarrow W_{n-1}\longrightarrow
\ldots\longrightarrow W_1\longrightarrow W_0=W$ so that conditions (i) and
(ii)  Theorem \ref{MainTh} hold for the strict transform of $X$ at $W_n$.
\item[(2)] A pattern of manipulation of equations defining $X$, so as to
obtain, at least locally at an open covering of $W_n$, equations defining
the strict transform  $X_n$ of $X$  at $W_n$.
\end{enumerate}
 
So (2) indicates how the original equations defining $X$ have to be
treated at an affine open subset of $W_n$ in order to obtain  local equations
defining $X_n$. While this is very complicated in Hironaka's line of proof, 
here it is a  direct consequence of Theorem \ref{MainTh} (iii). In fact,
for algorithms that follow Hironaka's proof, to get   both  (1) and (2)
one  must consider the {\em strict transform} of the
ideal of the subscheme at each monoidal transformation. In that setting
one has to  choose  a {\em standard basis} of the ideal, which is a system of
generators of the ideal  of the subscheme, but such choice of generators
must be changed if the maximum Hilbert Samuel invariant drops  in the
sequence  of monoidal transformations. All of  these complications are  avoided in our new
proof,  which simplifies both (1) and (2). The simplifications attained  in (2) are illustrated 
by the following corollary of Theorem \ref{MainTh}:

\begin{corollary}
\label{Koszul} Under the assumptions and with the notation of Theorem \ref{MainTh},  
if $X$ is a complete intersection, then  the  resolution of ${\mathcal O}_X$ in terms of free ${\mathcal O}_W$-modules 
$$\ldots \longrightarrow \land^k{\mathcal O}_W\longrightarrow \land^{k-1}{\mathcal O}_W\longrightarrow \ldots, $$
induces the resolution of ${\mathcal O}_{W_r}$
$$\ldots \longrightarrow {\mathcal {L}}^{-k}\land^k{\mathcal O}_{W_r}\longrightarrow {\mathcal {L}}^{-k+1}\land^{k-1}{\mathcal O}_{W_r}
\longrightarrow \ldots,$$
in terms of locally free ${\mathcal O}_{W_r}-$modules.
\end{corollary}

\noindent {\em Proof: } Assume that $W$ is affine, and that $X\subset W$ is defined by
the complete
intersection ideal $ I(X)=<f_1,f_2,...,f_r> \subset \calo_W $. A resolution of
 $\calo_X$ by free $\calo_W-$modules can be defined in terms
of a Koszul complex. This complex is defined by taking the tensor product of
$$ C_i \ :=   \ 0 \longrightarrow  \calo_W.e_i \longrightarrow  \calo_W
\longrightarrow 0, $$
where each such complex $ C_i$ is defined by $e_i \longrightarrow f_i$,
for $i=1,...,r$.
 
In the setting of \ref{MainTh} we have that $\calo_W \subset \calo_{W_r}$ and that 
$f_i \in  {\mathcal {L}} \subset \calo_{W_r}$. In particular each $C_i$ induces a complex 
$$ \underline{C}_i \ :=   \ 0 \longrightarrow  {\mathcal {L}}^{-1}.e_i
\longrightarrow  \calo_{W_r}
\longrightarrow 0 $$
Finally note that Theorem \ref{MainTh} (iii) says that the tensor product of these
defines a resolution of $\calo_{X_r}$ in terms of locally free
$\calo_{W_r}$-modules. \qed
 
\

The curve of our example $C\subset
W=\ma_{\rac}^3 $ is a complete intersection; note that the result in this corollary
{\em will never} hold for a
desingularization of this curve given within Hironaka's line of proof, since
condition (iii) of our Theorem will fail  as seen in {\bf A)}.

\section{Basic objects}
\label{secbo} To achieve our results we   use   the notions of
{\em basic objects} and {\em resolution of basic objects} (Definitions \ref{DefBasic} and
\ref{defResol}, see also \cite{EncVil97:Tirol}). Both  Embedded
 Desingularization  and  Strong
Principalization of Ideals  can be obtained
 from a resolution of  suitably defined  basic objects.

\begin{definition}
\label{DefBasic}
    A \emph{basic object}  is a triple that consists of  a pair
    \( (W,E) \), an ideal
    \( J\subset\calo_{W} \) such that $(J)_{\xi}\neq 0$
for any \( \xi\in W \),  and a positive integer \( b \).
It is denoted by \( (W,(J,b),E) \). If the dimension of $W$ is $d$, then \( (W,(J,b),E) \) is said to
be a {\em $d-$dimensional basic object}. 
\end{definition}
\begin{definition}
{\rm The {\em singular locus} of
a  basic object is the closed set:
\[ \Sing(J,b)=\{\xi\in W\mid \nu_{J}(\xi)\geq b\}
    \subset W, \]
where $\nu_{J}(\xi)$ denotes the order of the ideal 
$J$ at the point $\xi$.} 
\end{definition}
\begin{definition}
 {\rm  A regular closed subscheme \(
Y\subset W \) is  \emph{permissible} for \( (W,(J,b),E)
\) if \( Y \) is permissible for the pair \( (W,E) \) and \(
Y\subset \Sing(J,b)\). Given $Y\subset \Sing(J,b)$ we
define a {\em transformation of basic objects} in the following
way: Let $E=\{H_1,\ldots,H_r\}$, and consider the blowing-up with center $Y$,  
$$W\longleftarrow W_{1}.$$  
Denote by \(H_{r+1}\subset W_{1} \) the exceptional hypersurface.
This blowing-up induces a transformation of pairs
$$(W,E)\longleftarrow (W_1,E_1)$$
as in Definition \ref{pairs}. If \( Y \) is irreducible and if \( c_{1}\) is the
order of $J$ at the generic point of \( Y\) (i.e. $\nu_J(Y)=c_1\geq
b$), then there is an ideal \( \overline{J}_{1}\subset\calo_{W_{1}}
\) such that
        \begin{equation}
\label{extra00}
 J\calo_{W_{1}}=I(H_{r+1})^{c_{1}}\overline{J}_{1}.
\end{equation}
The ideal $\overline{J}_1$ is usually called the {\em weak transform} of $J$. 
Note here that $\overline{J}_1$ does not vanish along $H_1$ (i.e.
$H_1\subsetneqq V(\overline{J}_1)$). Under these assumptions we define the ideal
\begin{equation}
\label{suerte1}
J_{1}=I(H_{1})^{c_{1}-b}\overline{J}_{1}
\end{equation}
and we set
$$(W,(J,b),E)\longleftarrow (W_{1},(J_{1},b),E_{1})$$
as the \emph{transformation of the basic object} $(W,(J,b),E)$}.
\end{definition}

In general, given a sequence of
transformations of basic objects 
\begin{equation}
\label{resol}
(W_{0},(J_{0},b),E_{0})
\longleftarrow\cdots\longleftarrow (W_{k},(J_{k},b),E_{k}), 
\end{equation}
at centers $Y_i$, $i=0,1,\ldots k-1$, we obtain expressions 
\begin{equation}
\label{aaa1}
 J_{i}=I(H_{r+1})^{a_{1}}\cdots I(H_{r+i})^{a_{i}}\overline{J}_{i}
\end{equation}
and 
\begin{equation}
\label{dosestrellas}
J_{0}{\mathcal O}_{W_i}=I(H_{r+1})^{c_{1}}\cdots
I(H_{r+i})^{c_{i}}\overline{J}_{i},
\end{equation}
with $c_j> a_j\geq 0$ (we refer to \cite[4.8]{EncVil97:Tirol} for a precise description of these exponents). 
 Thus we have the factorization 
\begin{equation}
\label{tresestrellas}
J_{0}{\mathcal O}_{W_i}={\mathcal M}_i J_i={\mathcal L}_{i}\overline{J}_{i}
\end{equation}
for suitable defined invertible ideals ${\mathcal M_i}$ and ${\mathcal L}_{i}$.
\begin{definition}
\label{defResol}
Sequence \ref{resol} is a \emph{resolution  of}
    \( (W_{0},(J_{0},b),E_{0}) \) if \( \Sing(J_{k},b)=\emptyset \).
\end{definition}
 
\begin{remark}
\label{so} If (\ref{resol}) is a resolution of the basic object $(W_0,(J_0,b),E_))$ then
$W_k\longrightarrow W_0$ defines an isomorphism over $W_0\setminus
V(J_0)$, and $J_0{\mathcal O}_{W_k}={\mathcal M}_kJ_k$, 
where ${\mathcal M}_k$ is an invertible sheaf of ideals, and $J_k$ has no points of
order $ \geq b$ in $W_k$. In particular a Strong Principalization of a sheaf of ideals
 $I\subset {\mathcal O}_{W}$ follows
by taking
a resolution of a  basic object $(W_0,(J_0,b),E_0)$  with $W_0=W$, $J_0=I$, $b=1$ and $E_0=\emptyset$.
\end{remark}

A resolution of basic objects is usually approached by means of an
{\em algorithm of resolution of basic objects}, and the factorization 
in (\ref{aaa1}) plays a central role in this procedure (see Section \ref{last1}). To prove Theorem \ref{MainTh} we modify the 
algorithm of resolution
of basic objects which appears  in \cite{EncVil98}.

\section{Proof of  Theorem \ref{MainTh}}
\label{sketch}
\begin{parrafo}\label{dos}{\rm Let $(W_0,E_0=\emptyset)$ and  $X$ be  as in
Theorem \ref{MainTh} and let $J_0=\cali(X)\subset {\mathcal O}_{W_0}$. To prove
Theorem \ref{MainTh} we will define a suitable sequence of transformations of
pairs at permissible centers
\begin{equation}\label{eq1}
(W_0,E_0=\emptyset)\longleftarrow \ldots\longleftarrow (W_m,E_m),
\end{equation}
together with expressions
\begin{equation}\label{factorization}
J_0{\mathcal O}_{W_i}=\cali(H_1)^{c_1}
\cdot\ldots\cdot\cali(H_i)^{c_i}\overline{J}_i
\end{equation}
 as in  (\ref{dosestrellas}), so that the conclusion of Theorem
\ref{MainTh} holds for
$$J_0{\mathcal O}_{W_m}= \cali(H_1)^{c_1}\cdot\ldots\cdot\cali(H_k)^{c_m}
\overline{J}_m.$$

\

First we will  motivate the  argument of the proof of Theorem \ref{MainTh} 
from the algebraic point of view. This will 
lead us to  Definition \ref{omni} and to the statement of
Lemma \ref{lema}. The proof of Theorem \ref{MainTh} is presented in
\ref{proofoftheo}.

\

The proof of Lemma \ref{lema} will be given in
Section \ref{last}. In order to prove  this lemma  we will need some extra
ingredients which  will be presented in Section \ref{last1}. In Section
\ref{examples} we exhibit some examples illustrating the main ideas  of
Section \ref{last1}.}
 
\end{parrafo}
 
\begin{parrafo}\label{algebraicview}{\rm
 {\bf Conditions on  (\ref{eq1}) and (\ref{factorization}): The
embedded relative codimension.} 
Let  $(W,E=\{H_1,..,H_s\})$ be a pair and let  $y\in W$ be a point. Then 
we can set a
regular system of parameters $\{x_1,\ldots,x_n\}\subset {\mathcal O}_{W,y}$, so
that if $y\in \cup_{j=i_1}^{i_r}H_{i_j}$, then $\cup_{j=i_1}^{i_r}H_{i_j}$ is
locally defined by
\begin{equation}
\label{bdd}x_{i_1}\cdot\ldots\cdot x_{i_r}=0
\end{equation}
at a neighborhood of $y$ in $W$.

\
 
Assume that (\ref{eq1}) is defined so that at any point $y\in
V(\overline{J}_m)\subset W_m$,
 $(\overline{J}_m)_y$ is an ideal of order
1 in the local regular ring ${\mathcal O}_{W_m,y}$. In such case we may
choose a regular system of parameters
$\{z_1,\ldots,z_n\}\subset {\mathcal O}_{W_m,y}$
so that 
\begin{equation}
\label{add}
z_1\subset (\overline{J}_m)_y.
\end{equation}

\
 
The proof of  Theorem \ref{MainTh} is based in the idea of 
defining  suitable sequences of transformation of basic 
objects 
\begin{equation}
\label{eq1prima}
(W_0,E_0=\emptyset)\longleftarrow \ldots\longleftarrow (W_k,E_k),
\end{equation}
such that
at any point
$y\in V(\overline{J}_k)\subset W_k$, there is a regular system
of parameters so that both conditions (\ref{bdd}) and (\ref{add}) are
simultaneously
satisfied. This leads to the following definition: }
\end{parrafo}

\begin{definition}
\label{omni} {\rm Let $(W,E)$ be a pair, and let
$I\subset {\mathcal O}_{W}$ be a non-zero sheaf of ideals.
 
\noindent
 
(1) We shall say that $I$ has {\em relative embedded local codimension }
$\geq a$ {\em at a point} $y\in W$, if either $I_y={{\mathcal O}}_{W,y}$, or
there is a regular system of parameters $\{x_1,x_2, \ldots , x_n \} $ at
${{\mathcal O}}_{W,y}$, such  that:
 
\begin{enumerate}
\item[(i)]
$<x_1,x_2, \ldots , x_a>\subset I_y \subset {{\mathcal O}}_{W,y}$,
and
 
\item[(ii)] any hypersurface $H_i\in E$ containing the point $y$ has a
local equation 
$$\cali(H_i)=<x_{i_j}> \subset {{\mathcal O}}_{W,y},$$
with $i_j > a$.
\end{enumerate}
 
(2) We shall say that $I$ has
{\em relative embedded codimension } $\geq a$ in $(W,E)$, if both conditions in (1)
hold at every  point $y \in W$. We will abbreviate this saying that $I$ has 
{\em  relative codimension } $\geq a$ in $(W,E)$.} 
 
\end{definition}
 
\begin{remark}\label{rem2} Note
that if $I$ is of relative codimension $\geq a$, then the closed
subscheme $V(I)$ is in fact of codimension $\geq a$ in $W$. Note also that any
non-zero ideal $I \subset {{\mathcal O}}_{W}$ is of relative codimension
$\geq 0$ in $(W,E)$ since such condition is empty.
\end{remark}

\begin{remark}
If $X\subset W$ is a regular scheme of pure codimension $e$, then the sheaf of ideals 
$\cali(X)\subset {\mathcal O}_{W}$ has relative codimension $\geq e$ in $(W,E=\emptyset)$. 
\end{remark}
 
\begin{remark}
\label{rm3} 
Let $(W_i,(J_i,b),E_i)$ be a basic object and let 
$J_i=\call_i\overline{J}_i$ be as in (\ref{aaa1}). 
 If $\overline{J}_i$ has relative
codimension $\geq a$ in $(W,E)$, then for every $y\in V(\overline{J}_i)$
there is a
regular system of parameters $\{x_1,\ldots,x_n\}\subset {\mathcal O}_{W_i,y}$
such that
\begin{equation}
\label{rm4}\begin{array}{c}J{\mathcal O}_{W_i,y}=x_{i_1}^{c_1}\cdot\ldots
\cdot x_{i_r}^{c_r}(\overline{J}_i)_y; \mbox{ }
\mbox{ } \left<x_1,\ldots,x_a\right>\subset (\overline{J}_i)_y
\end{array}
\end{equation}
with $a<i_1<i_2<\ldots<i_r\leq n$ and $(\call_i)_y=<x_{i_1}^{c_1}\cdot\ldots
\cdot x_{i_r}^{a_r}>$.
\end{remark}
 
\begin{definition}
\label{bocodim} With the same  notation as in  Remark \ref{rm3}, 
we will say that $\overline{J}_i$ is of $(W_i,E_i)-$ {\em codimension
$\geq a$} if $\overline{J}_i$ has relative codimension $\geq a$ in $(W_i,E_i)$.
\end{definition}
 
Now assume that  $X$ is under the
assumptions  of  Theorem \ref{MainTh}, and that it has codimension $e$ as
subscheme in $W$.
Then the local codimension of $J_0=I(X)$  is
$\geq e$ at every 
 point of $V_0=W_0\setminus\Sing(X)$.  Therefore  $J_0$ will be of $(W_0,E_0)-$codimension $\geq a$ 
for some $a \leq e$ (since by Remark \ref{rem2} we may always assume $a=0$). Now
Theorem \ref{MainTh} will follow from Lemma \ref{lema}
and the argument given  in \ref{proofoftheo}.
\begin{lemma}
\label{lema}
Let $X$ be under the assumptions of Theorem \ref{MainTh}, let
$$(W_0, (J_0,1),E_0)=(W,(\cali(X),1),\emptyset)$$  
and  assume that there is a
sequence of transformations of pairs 
\begin{equation}
\label{2eq2}(W_0,E_0)\longleftarrow  \ldots\longleftarrow
(W_k,E_k)
\end{equation}
such that 
$$\cali(X){\mathcal O}_{W_k}=\call_{k}
\overline{J}_{k}$$ 
and that  $\overline{J}_{k}$ is of  $(W_k,E_k)-$codimension
$\geq a$. Let $ U_k \subset W_k$ be a non-empty open set such that
for every $y \in U_k$, the local codimension of $(\overline{J}_{k})_y $ is
$ \geq a + 1$. Then we can define an enlargement of the sequence \ref{2eq2},
 \begin{equation}
\label{eq2000}
(W_k,E_k)\longleftarrow
\ldots\longleftarrow (W_N,E_N)
\end{equation}
so that:
\begin{enumerate}
\item[(i)] $\cali(X){\mathcal O}_{W_N}=\call_{N}\overline{J}_{N}$and
$\overline{J}_{N}$ is of $(W_N,E_N)-$codimension $ \geq a+1$.
\item[(ii)] The birational morphism $W_k \longleftarrow W_N$ defines
an isomorphism over the open set $U_k\subset W_k$.
\end{enumerate}
\end{lemma}

\begin{parrafo}
\label{proofoftheo}
 {\rm {\bf  Proof of  Theorem \ref{MainTh}.}  Let $X\subset W$ be a
subscheme of codimension e, under the
assumptions of Theorem \ref{MainTh}. Set $W_0=W$, $J_0=\cali(X)$,
$E_0=\emptyset$, and $U_0=W_0 \setminus \Sing (X)$. Then $U_0\subset W$ is a 
dense open  subset, $U_0\cap X$ is dense in $X$, and $X$ is regular at all 
points of
$U_0\cap X$. Since $E_0= \emptyset$, 
the ideal $J_0$ is of relative local codimension
$\geq e$ at at every point $x\in U_0\cap X$. In particular, there is an integer $a$, $ 0
\leq a \leq e$,  so that $J_0$ has
$(W_0,E_0)$-codimension $\geq e$ (see Remark \ref{rem2}). If $a < e$, by successive applications of 
Lemma \ref{lema}, we may define a  sequence of
transformations
 \begin{equation}
\label{transformation}
(W_0,E_0)\longleftarrow
\ldots\longleftarrow (W_{N_e},E_{N_e})
\end{equation}
so that if 
\begin{equation}
\label{jotae}
 I(X){\mathcal O}_{W_{N_e}}={\call}_{N_e}\overline{J}_{N_e}
 \end{equation}
then 
\begin{enumerate}
\item[{(I)}] $\overline{J}_{N_e}$ has
$(W_{N_e},E_{N_e})-$codimension $\geq e$.
 \item[{(II)}] $W_0\longleftarrow W_{N_e}$ defines an isomorphism over $U_0$, say
$U_0\simeq U_{N_e}\subset W_{N_e}$.
\end{enumerate}

\
  
First we will make some elementary geometric remarks on the closed set
$V(\overline{J}_{N_e})\subset W_{N_e}$:
 
\begin{enumerate}
\item[(a)]
$V(\overline{J}_{{N_e}})$ has codimension $\geq e$ as closed subscheme of
$W_{{N_e}}$. Let
\begin{equation}
\label{compo}
V(\overline{J}_{{N_e}})=F_1\cup \ldots \cup F_l\cup C_{l+1}\cup
\ldots\cup C_p
\end{equation} 
be the union of irreducible components (each of codimension at
least $e$ in $W_{N_e}$), where the $F_i$ are the components of codimension
$e$. Set $F=F_1\cup \ldots\cup F_l$.

\item[(b)] Each component $F_i$  is a
smooth and {\em connected component} of
$V(\overline{J}_{N_e}{})$.
 
\item[(c)] $V(\overline{J}_{N_e}{})
\cap U_{N_e} \simeq X\setminus \Sing(X)=X \cap U_0$ ($U_{N_e}$ and $U_0$ as
in (II)).
 
\end{enumerate}
 
Conditions (a) and (b) can be checked from Definition \ref{omni}. Condition
(c) follows from (II). If we assume that $X$ is of
{\em pure} codimension $e$, then (b) and Definition \ref{omni} assert that:

\begin{enumerate}
\item[(d)] If $X_{N_e}\subset W_{N_e}$ denotes the strict
transform of $X$, then 
$$X_{N_e}=F_1\cup\ldots\cup F_{l^{\prime}}, \mbox{ } \mbox{ }l^{\prime}\leq l$$ 
is a disjoint union of closed regular subschemes.
\end{enumerate}
Condition (ii)  in Definition \ref{omni} asserts that  conditions (i)
and (ii) of
Theorem \ref{MainTh} hold for $W_{N_e}\longrightarrow W_0$. Now by (d)
we conclude that:
 
\begin{enumerate}
\item[(e)] For every 
$y\in X_{N_e}=F_1\cup\ldots\cup F_{l^{\prime}}$,
$$(\overline{J}_{{N_e}})_y=(\cali(X_{N_e}))_y\subset {\mathcal O}_{W_{N_e},y}.$$
\end{enumerate}
In fact, we know that $I(X_{N_e})_y$ is a primary component of
$(\overline{J}_{{N_e}})_y$, so  $(\overline{J}_{{N_e}})_y\subset
(\cali(X_{N_e}))_y$; if this inclusion were proper, then
$V(\overline{J}_{{N_e}})$ would have codimension $> e$ locally at $y$
(but $y\in F\subset V(J_{N_e})$). Note that (e) is saying that condition
(iii) of Theorem \ref{MainTh} holds locally at each point in
$F_1\cup\ldots\cup F_{l^{\prime}}$.
 
\
 
Since $V(\overline{J}_{{N_e}})$ is a {\em disjoint union } of
$F_1\cup \ldots \cup F_{l^{\prime}}$ and
$F_{l^{\prime}+1}\cup \ldots \cup F_l\cup C_1\cup \ldots \cup C_p$, we can
express $\overline{J}_{{N_e}}$ as a product of two ideals,
$\overline{J}_{{N_e}}=[A_{N_e}]_1\cdot[A_{N_e}]_2$,  so that
$$V([A_{N_e}]_1)=F_1\cup \ldots \cup F_{l^{\prime}}
\mbox{ and }V([A_{N_e}]_2)=F_{l^{\prime}+1}\cup \ldots \cup F_l\cup C_1\cup
\ldots \cup C_p.$$
Note finally that $V([A_{N_e}]_2)\subset W_{N_e}\setminus V([A_{N_e}]_1)$.
If $X$
is of pure codimension $e$, then (iii) of Theorem \ref{MainTh} is finally
achieved by setting$$(W_{N_e},E_{N_e})\longleftarrow \ldots\longleftarrow
(W_r,E_r)$$ so as to define a principalization of
$[A_{N_e}]_2\subset {\mathcal O}_{W_{N_e}}$.

\
 
If $X$ is  non-pure-dimensional the proof 
follows in a rather similar fashion: By (b) we have that at $W_{N_e}\setminus V([A_{N_e}]_1)$  the ideal 
$[A_{N_e}]_2$ factors as a product of two sheaves of ideals
\begin{equation}
\label{over}
[A_{N_e}]_2=[A_{N_e}]_3\cdot [A_{N_e}]_4,
\end{equation}
so that 
\begin{equation}
\label{over1}
V([A_{N_e}]_4)=C_{l+1}\cup\ldots\cup C_p
\end{equation}
where $\mbox{dim}(C_i)\geq e+1$, for $i=l+1,\ldots,p$. 

\

Now $X_{N_e}\cap \left[W_{N_e}\setminus V([A_{N_e}]_1)\right]$ is a union of some of the 
irreducible components 
$C_i$ of $V([A_{N_e}]_4)$ appearing in (\ref{over1}). These components are regular in a dense 
open set. Then we apply Lemma \ref{lema}  for some $e^{\prime} > e$ and make use
of (b) to separate those $e^{\prime}$-dimensional components of
$V([A_{N_e}]_4)$ which are components of the strict transform of X from those
which are not. Thus 
 we repeat the same argument that we applied in the pure codimensional case, but now for higher 
codimensions.
\qed}\end{parrafo}

\section{Algorithms of resolution of basic objects}
\label{last1}
Lemma \ref{lema} will follow from a suitable Algorithm of
Resolution of Basic Objects.
An  Algorithm of Resolution of Basic Objects is defined by means of  an
{\em upper semi-continuous function}
(see Definition \ref{AlgResol}). 
In \ref{secinduction} we will briefly describe  these upper semi-continuous
functions to give an idea of the strategy that we will 
follow (see also \cite{EncVil97:Tirol}). 
The novelty of our development appears in \ref{svive} and 
\ref{nuestroalgo}, where we show how the two main invariants
involved in the algorithm of resolution of basic objects (namely in the definition
 of the upper semi-continuous function) relate to our notion of relative
codimension stated in Definition \ref{omni}.

\begin{definition} 
\label{semicon}
let $X$ be a topological space, and let $(T, \geq )$ be 
 a totally ordered set. Let $g: X \longrightarrow T$ be an 
upper semi-continuous function, and assume that 
$g$ takes only finitely many values. 
The largest value achieved by $g$ will be denoted by 
    $$ \max g.$$
Clearly the set 
$$ \Max g=\{ x\in X : g(x)= \max g \} $$
is a closed subset of $X$.
\end{definition}

\begin{definition} 
\label{AlgResol}
 Let \( d \) be a non-negative integer. An {\em algorithm of resolution for
\( d \)-dimensional basic objects},  consists of:
 
\begin{enumerate}
 
\item[{\bf A.}]  A totally ordered set \( (I_{d},\leq) \).

\item[{\bf B.}]  For each $d-$dimensional basic object 
\( (W _ {0},(J_{0},b),E_{0}) \) (with $E_0$ is not necessarily empty),  an upper semi-continuous
 function 
\[ f_{0}^{d}:\Sing(J_{0},b)\longrightarrow I_{d} \] 
such that
\( \Max{f_{0}^{d}}\)  is permissible for \( (W_{0},(J_{0},b),E_{0}) \).
\item[{\bf C.}] Suppose that we define a sequence of 
  blowing-ups at permissible centers \(
Y_{i
}\subset\Sing(J_{i},b) \), \( i=0,\ldots,r-1 \), 
  \begin{equation}
\label{eqMainSeqResol}
  (W_{0},(J_{0},b),E_{0})\longleftarrow\cdots
\longleftarrow
 (W_{r-1},(J_{r-1},b),E_{r-1})
\longleftarrow
(W_{r},(J_{r},b),E_{r}),
\end{equation}
together with a sequence of  upper semi-continuous functions
$$f_{i}^{d}:\Sing(J_{i},b)\longrightarrow I_{d}, \mbox{ } 
i=0,\ldots,r-1,$$ 
so that 
$$Y_{i}=\Max{f_{i}^{d}}.$$
Then  if \( \Sing(J_{r},b)\neq\emptyset\),  there is an upper semi-continuous 
 function \[ f_{r}^{d}:\Sing(J_{r},b)\longrightarrow I_{d} \]
such that \( \Max{f_{r}^{d}} \) is permissible for
\( (W_{r},(J_{r},b),E_{r}) \).
\item[{\bf D.}]  For some index \( r_0 \),
 depending on the basic object \( (W_{0},(J_{0},b),E_{0}) \),
 the  sequence constructed in  {\bf C} is a
 resolution (i.e. \( \Sing(J_{r_0},b)=\emptyset \)).
\item[{\bf E.}]{\bf Properties.} The following properties should hold:
 \begin{enumerate}
 \item[i.]  If \( \xi\in \Sing(J_{i},b) \) and  \( \xi\not\in Y_{i} \)
  for
  \( i=0,\ldots,r_0-1 \),
  then
\( f_{i}^{d}(\xi)=f_{i+1}^{d}(\xi') \)
  via the natural identification of the point \( \xi \) with a
  point \( \xi' \) of \( \Sing(J_{i+1},b) \).
\item[ii.]  The resolution is obtained by transformations with
    centers \( \Max{f_{i}^{d}} \), for \( i=0,\ldots,r_0-1 \), and
    \[ \max{f_{0}^{d}}>\max{f_{1}^{d}}>\cdots >\max{f_{r_0-1}^{d}} \]
\item[iii.]  If \( J_{0} \) is the ideal of a
    regular pure dimensional subscheme \( X_{0} \), if 
    \( E_{0}=\emptyset \),  and if \( b=1 \), then the
 function \( f_{0}^{d} \) is constant.
\item[iv.]  For every \( i=0,\ldots,r_0-1 \), the closed set
    \( \Max{f_{i}^{d}} \) is smooth, equidimensional, and its dimension
    is determined by the value \( \max{f_{i}^{d}} \).
\end{enumerate}
    \end{enumerate}
\end{definition}
\begin{remark}
\label{nose}
 Note that {\bf B} asserts that the setting of
(\ref{eqMainSeqResol}) holds for \( r=1 \), whereas {\bf C} says
that whenever \(\Sing(J_{r},b)\neq\emptyset \) there is an
enlargement of (\ref{eqMainSeqResol}) with center \(
Y_{r}=\Max{f_{r}^{d}}\).
Property {\bf E} (i) ensures that the algorithm commutes with open
restrictions: If we restrict to a non empty open subset of $W$, 
then the restriction of the algorithm gives the resolution of the
restriction of the basic object to the open subset.
 Property E (ii) says that the heart of the matter in the algorithm
is to find a good upper semi-continuous function that controls
the singular locus of the basic object, determines the center to blow-up
(see {\bf B}, {\bf
C} and {\bf E} (iv)) and guarantees that at each step the
{\em singularities} are getting better (property {\bf E} (ii)).
\end{remark}

\begin{parrafo}
\label{existence}{\rm An algorithm of resolution of basic
objects, as defined in \ref{AlgResol} and with the properties stated in
{\bf E}, is presented in  \cite[7.15]{EncVil97:Tirol}.
There are several ways in which embedded desingularization can be
achieved in terms of a fixed algorithm of resolution
 of
basic objects. It is indicated in \cite{EncVil97:Tirol} (see the addendum) how a
particular theorem of embedded desingularization can be easily
deduced as a byproduct from the properties listed in {\bf E}. 
In what follows we will briefly explain how the functions $f_i^d$ are defined, 
and for more details the reader is referred to \cite{EncVil97:Tirol}. }
\end{parrafo}

\begin{center}
{\bf How to define  the functions $f_i^d$.}
\end{center}
 
\begin{definition} \label{functiont} Let
$(W_0,(J_0,b),E_0=\{H_1,\ldots,H_l\})$ be a $d-$dimensional basic object.
Given a   sequence of permissible  transformations,
\begin{equation}
\label{extra11} 
(W_{0},(J_{0},b),E_{0})
\stackrel{\pi_1}{\longleftarrow}\cdots\stackrel{\pi_r}{\longleftarrow} (W_{r},(J_{r},b),E_{r})
\end{equation}
together with the corresponding expressions
\begin{equation}
\label{extra22} J_{i}=I(H_{l+1})^{a_{1}}\cdots
I(H_{l+i})^{a_{r}}\overline{J}_{i},
\end{equation}
where $H_{l+i}$ denotes the exceptional divisor at that $i-$th blowing-up,
we define the upper semi-continuous functions  $$\begin{array}{cccc}
\wordid: & \Sing(J_{i},b) & \longrightarrow &
        \frac{1}{b}\mathbb{Z}\subset\mathbb{Q}\\ & \xi & \longrightarrow &
        \frac{\nu_{\overline{J}_{i}}(\xi)}{b}    \end{array}$$
\end{definition}
\begin{remark} Note that for any $k\in {\mathbb N}_{\geq 1}$ a resolution of 
$(W_0,(J_0,b),E_0)$ is {\em equivalent} to a resolution of $(W_0,(J_0^k,kb),E_0)$. In fact it suffices to 
take $k-$th powers in the factorization in (\ref{extra22}). So it is reasonable to expect that the function 
$\word_i^d$ be compatible with this equivalence. This explains why the 
function $\word_i^d$ is defined taking the quotient by $b$ in the previous definition.
\end{remark}
\begin{remark}\label{centroper}
Note that  since \(
    \wordid \) is  upper semi-continuous  
the set
$$\Max\word_i^d=\{\xi\in \Sing(J_i,b):\word_i^d(\xi)=\mbox{max
}\word\}$$
is closed in $\Sing(J_{i},b)\subset W_i$. Note also that 
$\max \word_i^d=0$ if and only if $\overline{J}_i={\mathcal O}_{W_i}$. 
\end{remark}
\begin{remark}\label{centroper1}
We will only  consider sequences of transformations of
basic objects at permissible centers $Y_i$ such  as the ones defined in
(\ref{extra11}), with the additional constraint that
        \begin{equation}
        \label{adcons}
Y_{i}\subset\Max\word_{i}^d\subset\Sing(J_{i},b)
        \end{equation}
for $i=0,1,\ldots,r$. In such case
    \begin{equation} \label{IneqwordPt}
        \word_{i-1}^d(\pi_{i}(\xi_{i}))\geq\word_{i}^d(\xi_{i})
    \end{equation}
    for every \( \xi\in\Sing(J_{i},b) \),
 and  equality holds if $\xi\notin Y_i$.
 In particular    
\begin{equation} 
\label{Ineqword}
        \max\word_{0}^d\geq\cdots\geq\max\word_{r}^d
    \end{equation}
(see  \cite[4.12]{EncVil97:Tirol}).
\end{remark}

\begin{definition} \label{Deft}
    {\rm Consider a sequence of transformations of basic objects as in
    (\ref{extra11}), with $$Y_i\subset\Max\word_i^d$$ for
    $i=0,1,\ldots,r$. Pick $k\in \{0,1,\ldots,r\}$.
    If \( \max\word_{0}^d>0 \)
    let \( k_{0} \) be the smallest index so that
    $$\max\word_{k_{0}-1}^d>\max\word_{k_{0}}^d=\max\word_{k}^d$$
    (\( k_{0}=0
    \) if \( \max\word_{0}^d=\cdots=\max\word_{r}^d \)). 
Write  
$$E_{k}=E_k^{+}\sqcup E_{k}^{-}$$
where \(E_{k}^{-} \) is the set hypersurfaces of 
    \(E_{k} \)  which are strict transforms of hypersurfaces of \(
E_{k_{0}}
    \). Now  define
        \[ n_{k}^d(\xi)=\left\{
    \begin{array}{lll}
        \#\{H\in E_{k}\mid \xi\in H\} & {\rm if} &
\word_{k}^d(\xi)<\max\word_{k}^d \\
        \#\{H\in E_{k}^{-}\mid \xi\in H\} & {\rm if} &
\word_{k}^d(\xi)=\max\word_{k}^d.
    \end{array}
    \right. \]
    }
\end{definition}
 
\begin{definition}
\label{Deft1} {\rm With the hypothesis of   Definition
\ref{Deft}, if  \( \max\word_{r}^d>0 \), we define, for the index $r$, a
function $t_r^d$ by
setting:
\[
\begin{array}{ccccc}
 t_{r}^d: & \Sing(J_{r},b) & \longrightarrow &
    (\mathbb{Q}\times\mathbb{Z},\leq) &  \hbox{(lexicographic order)} \\
     & \xi & \longrightarrow & (\word_{r}^d(\xi),n_{r}^d(\xi)) &
\end{array}
\]
In the same way we define functions \(
t_{r-1}^d,t_{r-2}^d,\ldots,t_{0}^d \)}. If $Y_r\subset \Max t^d_r$
is a permissible center, then $Y_r$ is said to be a {\em
$t^d_r-$permissible center}.
\end{definition}

\begin{parrafo}
\label{PropiedFt}
{\rm {\bf Properties of the function \( t^d_i \)} 
(cf. \cite[4.15]{EncVil97:Tirol}). 
The functions  $$t_i^d:
\Sing(J_i,b)\longrightarrow I_d$$ have the following properties:
    \begin{enumerate}
        \item  Each \( t_{i}^d:\Sing(J_{i},b)\longrightarrow
        \mathbb{Q}\times\mathbb{Z} \) is upper semi-continuous.
        \item  If in the sequence of transformations 
\begin{equation}
\label{extra01000} 
\begin{array}{c}
(W_{0},(J_{0},b),E_{0})
\stackrel{\pi_1}{\longleftarrow}(W_{1},(J_{1},b),E_{1})\stackrel{\pi_2}
{\longleftarrow}\cdots\\
\cdots\stackrel{\pi_{r-1}}{\longleftarrow}(W_{r-1},
(J_{r-1},b),E_{r-1})\stackrel{\pi_r}
{\longleftarrow} (W_{r},(J_{r},b),E_{r})
\end{array}
\end{equation}
the permissible centers $Y_i$ have the
 additional constrain 
\[ Y_{i}\subset\Max{t_{i}^d}\subset\Max\word_{i}
\qquad i=0,1,\ldots,r-1\]
then for each index \( i=0,1,\ldots,r \),  
\[ t_{i-1}^d(\pi_{i}(\xi_{i}))\geq t_{i}^d(\xi_{i}) \qquad
      \mbox{ }   \mbox{ for all } \mbox{ } \xi_{i}\in\Sing(J_{i},b) \]
and equality holds if \( \pi_{i}(\xi_{i})\not\in Y_{i-1} \).
        In particular
        \[ \max{t_{0}^d}\geq\cdots\geq\max{t_{r}^d}. \]
\item  We say that \( \max{t}^d \) drops at \( i_{0} \) if \(
        \max{t_{i_{0}-1}^d}>\max{t_{i_{0}}} \).
        If \( \max\word_{0}^d=\dfrac{b'}{b} \) and
\( \dim{W_{0}}=d \), note
        that \( \max{t_{i}^d}=\left(\dfrac{s}{b},m\right) \),
\( 0\leq s\leq b'
        \), \( 0\leq m\leq d \).
        So it is clear that \( \max{t^d} \) can drop at most
\( b'd \) times.
\item The functions $t^d_i$ are the {\em inductive invariant}:
via some form of induction, which will be described in \ref{secinduction},      
it is possible to  construct a unique enlargement of
        (\ref{extra01000}) 
        \begin{multline}
        \label{nueva2}
            (W_{0},(J_{0},b),E_{0})\longleftarrow\cdots\longleftarrow
            (W_{r},(J_{r},b),E_{r})\longleftarrow\hspace{2cm} \\
            \longleftarrow
            (W_{r+1},(J_{r+1},b),E_{r+1})\longleftarrow\cdots\longleftarrow
            (W_{N},(J_{N},b),E_{N})
        \end{multline}
        so that
        \( \max{t_{r}^d}=\max{t_{r+1}^d}=\cdots=\max{t_{N-1}^d} \)
        and either
        \begin{enumerate}
            \item  \( \Sing(J_{N},b)=\emptyset \).
 
            \item  \( \Sing(J_{N},b)\neq\emptyset \)
and \( \max\word_{N}=0 \).
 
            \item  \( \Sing(J_{N},b)\neq\emptyset \),
\( \max\word_{N}>0 \)
            and \( \max{t_{N-1}}>\max{t_{N}} \).
        \end{enumerate}
\end{enumerate}
 
Note that Property 3 says that $\max\word^d$ can drop finitely many
times, in particular for some index $N$, either (a) or (b) will
hold.}
\end{parrafo}

\begin{remark}
 \label{algomas} We will consider here
$d-$dimensional basic objects $(W_0,(J_0,b),E_0)$ together with sequences of
$t_i^d-$permissible transformations,
\begin{equation}
\label{algoextra} (W_{0},(J_{0},b),E_{0})
\longleftarrow\cdots\longleftarrow (W_{r},(J_{r},b),E_{r}).
\end{equation}
In such case $\max \mbox{w-ord}_i^d\geq \max \mbox{w-ord}_{i+1}^d$ and
$\max t_i^d\geq \max t_{i+1}^d$. We shall also consider the factorization
\begin{equation}
\label{masextra}J_i=\cali(H_{l+1})^{a_1}\cdot\ldots\cdot\cali(H_{l+i})^{a_i}
\overline{J}_i\end{equation}as in (\ref{aaa1}).
 
The condition
$\max \mbox{w-ord}_i^d>0$ is equivalent to $\overline{J}_i\subset
{\mathcal O}_{W_i}$, and  $\overline{J}_i\neq {\mathcal O}_{W_i}$. If
$\max \mbox{w-ord}_i^d=0$, then $\overline{J}_i={\mathcal O}_{W_i}$ and a
resolution of $(W_i,(J_i,b),E_i)$ is a simple combinatorial problem (cf. 
\cite[Section 5]{EncVil97:Tirol}). If
$\max \mbox{w-ord}_i^d>0$, then
\begin{equation}\label{ulyi}
\Max t_i^d\subset
\Max \mbox{w-ord}_i^d \subset V(\overline{J}_i).
\end{equation}
 \end{remark}
 
\begin{remark}
\label{claro} 

Let  $(W_0,(J_0,1),E_0=\{H_1,\ldots,H_l\})$ be a
$d-$dimensional basic object and consider a  sequence of permissible
transformations,
\begin{equation}
\label{extra0100} (W_{0},(J_{0},1),E_{0})
\longleftarrow\cdots\longleftarrow (W_{r},(J_{r},1),E_{r})
\end{equation}
together with the corresponding expressions, 
\begin{equation}
\label{2extra22} J_{i}=I(H_{l+1})^{a_{1}}\cdots
I(H_{l+i})^{a_{r}}\overline{J}_{i}, 
\end{equation}
then  if $\max t_i^d=(1,0)$, $\overline{J}_i$ is of $(W_i,E_i)-$codimension
$\geq 1$ (see Remark \ref{svive} for the proof of this fact). Therefore, in our development 
we are particularly interested in resolutions of basic objects $(W,(J,b)E)$ with $b=1$. However, 
the construction of such resolutions involves the definition of new basic objects in order to get 
 suitable invariants, and then we will have to consider the case when $b>1$ (see \ref{secinduction}).
\end{remark}

\begin{parrafo}
\label{secinduction}
{\rm {\bf The functions  $t^d_i$  and  the functions $f_i^d$.} 
We will present the construction of the functions $f_i^d$ in three steps, 
 and  we refer to \cite{EncVil97:Tirol} 
for more details. First 
we will give a brief idea of the strategy that we will follow: 
The main point is to give a precise formulation of the 
{\em inductive invariant} and the unique enlargement of  the sequence 
(\ref{nueva2}). 
To this purpose, 
 given a 
basic object  $(W_{i},(J_{i},b),E_{i})$  we will associate to it another 
$(W_i,(J''_{i},b''),E''_{i})$  which will be 
 constructed by means of a 
suitable differential operator applied to $J_i$. 
Using  this new basic object we can associate to $(W_{i},(J_{i},b),E_{i})$ 
a new $d-1$ dimensional basic object which will be described in step 2. 
This is the key point of induction on the dimension of the basic objects. 

\

Although the definition of these functions is based in a local argument
(see Steps 1 and 2) the construction is natural enough so as to produce
globally defined functions (see Remark \ref{versa}). 
 
\noindent {\bf Step 1.} Let $(W_i,(J_i,b),E_i)$ be a basic object, and let
$$t_i^d: \Sing(J_i,b)\longrightarrow ({\mathbb Q}\times {\mathbb Z},\leq ),$$
be as in Definition \ref{Deft1}. 
Consider the closed set $\Max t_i^d\subset\Sing(J_i,b)$. Then for
any $x\in \Max t_i^d$  there is a new basic object
$$(W_i,(J_i^{\prime \prime},b^{\prime \prime}),E_i^{\prime \prime})$$
(see (cf. \cite[9.5]{EncVil97:Tirol}, particularly (9.5.7)), with the following
properties:
\begin{enumerate}
\item[{\bf (a)}] $\Max{t_{i}^d}=\Sing(J''_{i},b'')$.
\item[{\bf (b)}] Any sequence of transformations at permissible centers
$C_{j}\subset\Sing(J''_{j},b'')$,       
\begin{equation}
        \label{nuevas44}
        (W_i,(J''_{i},b''),E''_{i})\longleftarrow\cdots\longleftarrow
(W_{s},(J''_{s},b''),E''_{s})
        \end{equation}
induces  a sequence  with same
centers $C_{j}\subset\Sing(J_{j},b)$,
        \begin{equation}
        \label{nuevas55}
        (W_{i},(J_{i},b),E_{i})\longleftarrow\cdots\longleftarrow
(W_{s},(J_{s},b),E_{s})
        \end{equation}
and has the following properties:
 
\begin{enumerate}
 
\item[(i)] If $\max{t_{i}^d}=\cdots=\max{t_{s-1}^d}$ then
$\Max{t_{j}^d}=\Sing(J''_{j},b'')$ for $i\leq j<s$.
 
\item [(ii)] $\Sing(J''_{i},b'')=\emptyset$ if and only
$\max{t_{s-1}^d}>\max{t_{s}^d}$.
\end{enumerate}
\end{enumerate}
\noindent {\bf Step 2.} Let $R(1)(\Max t_i^d)$ be the union of components of
$\Max t^d$ which have codimension 1 in $W_i$.
 
\begin{enumerate}
\item[{\bf (1)}] If $R(1)(\Max t_i^d)\neq \emptyset$ then
$R(1)(\Max t_i^d)$ is regular, and it is open
and closed in $\Sing(J^{\prime \prime},b^{\prime \prime})$. Furthermore, it has
normal crossings with $E$. In this case this is our canonical choice of center and 
our resolution function 
$f_i^d$ will be defined so that $\Max f_i^d=R(1)(\Max t_i^d)$. Therefore we
blow-up at $R(1)(\Max t_i^d)$,
\begin{equation}
        \label{nuevas555}
        (W_{i},(J_{i},b),E_{i})\longleftarrow(W_{i+1},(J_{i+1},b),E_{i+1}).
        \end{equation}
and either $\max t_i^d>\max t_{i+1}^d$ or 
$R(1)\Max(t_{i+1}^d)=\emptyset$. 

\item[{\bf (2)}] If $R(1)(\Max t_i^d)=\emptyset$ then we associate to
$(W_i, (J^{\prime \prime}_i,b^{\prime \prime}), E^{\prime \prime}_i)$  a
$d-1-$dimensional basic
object,
\begin{equation}\label{otramas}
(Z_i,(C(J^{\prime \prime}_i),b^{\prime \prime}!),
\overline E^{\prime \prime}_i),
\end{equation}
where $Z_i\subset W_i$ is a smooth
subscheme of dimension $d-1$, $\overline E_i^{\prime \prime}$ is the
restriction of
$E_i^{\prime \prime}$ to $Z_i$, and $C(J^{\prime \prime}_i)$ is the {\em coefficient
ideal} of $J^{\prime \prime}$
(cf. \cite[9.3]{EncVil97:Tirol})).
\end{enumerate}
This $d-1$ dimensional basic object has the following
properties:
 
\begin{enumerate}
 
\item[{\bf (a)}] $\Max{t_{0}^i}=
\Sing(C(J_i^{\prime \prime}),b^{\prime \prime}!)$.
\item[{\bf (b)}]
Any sequence of transformations        \begin{equation}
        \label{nuevas4}
        (Z_{i},(C(J_{i}^{\prime \prime}),b^{\prime \prime}!),\overline{E}_{i})
\longleftarrow\cdots
\longleftarrow
(Z_{s},(C(J_{s}^{\prime \prime}),b^{\prime \prime}!),\overline{E}_{s})  
\end{equation}
at permissible centers $D_j\subset
\Sing((C(J_{j}^{\prime \prime}),b^{\prime \prime}!)$
induces a  sequence
with same
centers $D_{j}\subset\Sing(J_{j},b)$,
        \begin{equation}
        \label{nuevas5}
        (W_{i},(J_{i},b),E_{i})\longleftarrow\cdots\longleftarrow
(W_{s},(J_{s},b),E_{s})
        \end{equation}
and
\begin{enumerate}
\item[(i)] If $\max{t_{i}^d}=\cdots=\max{t_{s-1}^d}$ then
$\Max{t_{j}^d}=\Sing(C(J_{j}^{\prime \prime}),b^{\prime \prime}!)$ for
$i\leq j<s$.
 
\item [(ii)] $\Sing(C(J_{i}^{\prime \prime}),b^{\prime \prime}!)=\emptyset$
if and only
$\max{t_{s-1}^d}>\max{t_{s}^d}$.
\end{enumerate}
\end{enumerate}
 
\noindent{\bf Step 3.} If $\max \mbox{w-ord}_i^j>0$ then 
Apply the whole process from Step 1, now to the
$d-1-$dimensional basic object
$(Z_{i},(C(J_{i}^{\prime \prime}),b^{\prime \prime}!)=
(Z_{i},(J_{i}^{d-1},b^{\prime \prime}!)$ so as
to define $$t_i^{d-1}:
\Sing(C(J^{\prime \prime}_i),b^{\prime \prime}!) \longrightarrow \tqz.$$

The successive application of steps 1, 2 and 3 leads us to define functions 
$$t_i^j: \Sing(J_i^j,b)\longrightarrow \tqz,$$
for $j=d,d-1,\ldots,1$. 
However 
note that the functions 
$$t_i^j:\Sing(J_i,b) \longrightarrow 
\tqz$$ 
are only defined if $\max \mbox{w-ord}_i^j>0$ (see Definition \ref{Deft1}). 
If \( \max\word_{i}^j=0
\) then
a resolution of the corresponding  basic object $(W_i,(J_i^j,b),E_i)$ 
is simple, and can also be achieved by
blowing up at the maximum of a suitable defined 
 function which values in a totally ordered set $(T,<)$ containing $\tqz)$, but this does not require 
a form of induction on the dimension. If we still
denote such function by $t_i^d$, say
 $$t^{j}_i:
\Sing(J_i^j,b) \longrightarrow (T,\leq),$$
then
the successive
application of Steps 1, 2 and 3  permits to define a function:
$$ f_i^d: \Sing(J_0,b)  \longrightarrow
\overbrace{T\times T\times\cdots\times T}^{k-\mbox{times}}, \mbox{ } \mbox{ }1\leq k\leq d$$
which is upper semi-continuous and verifies the conditions required in
Definition \ref{AlgResol}. We refer here to  in \cite[7.11]{EncVil97:Tirol} for a
precise description of this invariant. }
\end{parrafo}
 
\begin{remark}
\label{versa}
A major property of the
functions $t^d_{i}$ is that these
locally defined functions $t^{d-1}_0$ patch so as to define a
global function $t^{d-1}_i$, and the locally defined resolutions
(\ref{nuevas44}), (\ref{nuevas55}),  (\ref{nuevas4}) and (\ref{nuevas5}) are
sufficiently natural so as to define a (global) resolution of the basic object 
$(W_i,(J_i,b)E_i)$. This is discussed in \cite[9.5]{EncVil97:Tirol}.
\end{remark}

\begin{remark} \label{svive}
The proof of Lemma \ref{lema} is based on the notion of relative
codimension stated in Definitions \ref{omni} and  \ref{bocodim}. As we mentioned in 
Remark \ref{claro} this notion is  related to certain values of 
 the functions $t_i^d$:  
Let $X\subset W$ be under the assumptions of Theorem \ref{MainTh}. Then $X$ is
regular  in $X\setminus \Sing(X)$, which is a dense open set of $X$.
Consider
the basic object 
$$(W_0,(J_0,1),E_0)=(W,(\cali(X),1),\emptyset).$$ 
Note that the order of
$J_0$ is one at  every point $y\in X\setminus \Sing(X)$, and in fact  for every 
$y\in X\setminus \Sing(X)$,
$$t_0^d(y)=(1,0).$$ 
Note  also that $J_0$ is of
 local codimension $\geq 1$, for every 
$y\in X\setminus \Sing(X)$.
Consider the constructive resolution of the basic object
$(W_0,(J_0,1),E_0)$,
\begin{equation}
\label{gracias}
(W_0,(J_0,1)E_0)\longleftarrow \ldots\longleftarrow (W_N,(J_N,1),E_N).
\end{equation}
Since $\Sing{(J_N,1)}=\emptyset$, the ideal  $J_0{\mathcal O}_{W_N}$ defines a
principalization of $J_0$. Now, by property
{\bf E} of Definition \ref{AlgResol}  and since $\mbox{w-ord}_0^d(y)=1$ for
$y\in  X\setminus \Sing(X)$, it follows that there must be an index
$i_0\in \{0,1,\ldots,N\}$ in the sequence (\ref{gracias}) such that
$$\max \mbox{w-ord}_{i_0-1}>\max \mbox{w-ord}_{i_0}=1.$$
For such index $i_0$,
$$J_0{\mathcal O}_{W_{i_0}}=\call_{i_0}\overline{J}_{i_0}$$
and the order of $\overline{J}_{i_0}$ at every point
$y\in V(\overline{J}_{i_0})$ is one. Since the sequence (\ref{gracias})
is obtained by blowing-up at $t_i^d-$permissible centers, we may also choose
$i_1\geq i_0$ so that
     $$\max t_{i_1-1}^d>\max t_{i_1}^d=(1,0).$$
For such index $i_1$
$$
J{\mathcal O}_{W_{i_1}}=\call_{i_1}\overline{J}_{i_1}.
$$
 
We claim now that $\overline{J}_{i_1}$ has $(W_{i_1},E_{i_1})-$codimension
$\geq 1$. To see this we fix a point $y_{i_1}\in V(\overline{J}_{i_1})$. We 
will show  that the local condition of Definition \ref{omni} holds at
${\mathcal O}_{W_{i_1},y_{i_1}}$. Let $y_{i_0}$ denote the image of $y_{i_1}$
at $W_{i_0}$. By Remark \ref{centroper1}
$$
\mbox{w-ord}_{i_0}^d(y_{i_0}) \geq \mbox{w-ord}_{i_1}^d(y_{i_1}).
$$
Since 
$$
1=\max \mbox{w-ord}_{i_0}^d \geq \mbox{w-ord}_{i_0}^d(y_{i_0}) \geq
\mbox{w-ord}_{i_1}^d(y_{i_1})\geq 1,
$$
it follows that
$$
\mbox{w-ord}_{i_0}^d(y_{i_0})=\mbox{w-ord}_{i_1}^d(y_{i_1})=1.
$$
Since the order of $\overline{J}_{i_0}$ is one, 
we can choose a regular system of parameters
$$\{x_1,\ldots,x_n\}\subset {\mathcal O}_{W_{i_0},y_{i_0}}$$
so that $x_1\in \overline{J}_{i_0}$.
We still do not know if condition (ii) of (1) in Definition \ref{omni}
holds, since the hypersurfaces $H_j\in E_{i_0}$ may not be transversal to
$V(<x_1>)$.
Locally at $y_{i_0}$ set $Z_{i_0}=V(x_1)$; so $Z_{i_0}$
is a smooth hypersurface which contains $V(\overline{J}_{i_0})$. Let
$C_{i_0}$ denote the center of the transformation
$$
W_{i_0}\longleftarrow W_{i_0+1}.
$$
 
If $y_{i_0}\in C_{i_0}$ then
$$C_{i_0}\subset \Max t_{i_0}^d\subset V(\overline{J}_{i_0})
\subset Z_{i_0}$$
and hence, if $Z_{i_0+1}$ denotes the strict transform of
$Z_{i_0}$ in $W_{i_0+1}$, the exceptional divisor intersects $Z_{i_0+1}$
transversally. On the other hand, the law of transformation defining
$$
(W_{i_0},(J_{i_0},1),E_{i_0})\longleftarrow
(W_{i_{0}+1},(J_{i_{0}+1},1), E_{i_{0}+1})
$$
is such that $\cali(Z_{i_0+1})\subset
\overline{J}_{i_0+1}$. In particular the smooth hypersurface $\cali(Z_{i_0+1})$
is in the same setting as $Z_{i_0}$, namely,
$$
V(\overline{J}_{i_0+1})\subset Z_{i_0+1}.
$$
 
Since $i_0\leq i_1$, we repeat
this argument $(i_0-i_1)-$times, 
\begin{equation}
\label{jnumero}
(W_{i_0},E_{i_0})\longleftarrow (W_{i_0+1},
E_{i_0+1})\longleftarrow \ldots\longleftarrow (W_{i_1},E_{i_1})
\end{equation}
and at the end we get that for each index $ i_0 \leq j \leq i_1 -1 $
\begin{equation}
\label{talvez}
 C_j \subset Z_j, 
\end{equation}
where $C_j$ denotes the center of the monoidal transformation
$ W_j \longleftarrow W_{j+1}$ and  $Z_j \subset W_j$ denotes
the strict transform of $Z_{i_0}\subset W_{i_0}$. Also
$$
V(\overline{J}_{l})\subset Z_{l},
$$
for each index $ i_0 \leq l \leq i_1$. 
 Now $y_{i_1}\in V(\overline{J}_{i_1})\subset Z_{i_1}$, and
$$t_{i_1}^d(y_{i_1})=(\mbox{w-ord}_{i_1}^d(y_{i_1}), \mbox{n}_{i_1}^d(y_1))=
(1,0).$$
Thus  condition
$n_{i_1}^d(y_{i_1})=0$ asserts that none of the strict transforms of
hypersurfaces
$H_j\in E_{i_0}$ contains the point $y_{i_1}$, while all the exceptional
divisors in sequence (\ref{jnumero}) have normal crossings with
$Z_{i_1}$ by condition (\ref{talvez}).  Therefore $\overline{J}_{i_1}$ has
$(W_{i_1},E_{i_1})-$codimension $\geq 1$ at $y_{i_1}$. 
\end{remark}
 
\begin{remark}
\label{nuestroalgo}
{\bf Case $\max t^d=(1,0).$} With the same
notation as in
Remark \ref{svive}, set $i_1$ so that
$$\max t_{i_1-1}^d>\max t_{i_1}^d=(1,0).$$
Since $\max \mbox{w-ord}^d_{i_1}=1$, it follows that $\Max \mbox{w-ord}^d_{i_1}=
V(\overline{J}_{i_1})$. Furthermore, since $t^d_{i_1}$ is an
upper semi-continuous function, it also follows that
$\Max t^d_{i_1}=V(\overline{J}_{i_1})$. 
If $\max t_{i_1}^d=(1,0)$ then the basic object
$(W_{i_1},(J_{i_1}^{\prime \prime},b^{\prime \prime}),E_{i_1}^{\prime\prime})$
attached to
$\Max t^d_{i_1}$ as in (2) of Step 2 in  \ref{secinduction} is
defined by setting $b^{\prime \prime}=1$, 
$J_{i_1}^{\prime \prime}=\overline{J}_{i_1}$ and $E_{i_1}^{\prime\prime}=E_{i_1}$ 
(see  \cite[9.5]{EncVil97:Tirol}).

Write $J_{i_1}=\call_{i_1}\overline{J}_{i_1}$. Now if 
$y_{i_1}\in V(\overline{J}_{i_1})$ then there is a regular system of parameters
$\{x_1,\ldots,x_n\}\subset {\mathcal O}_{W_{i_1},y_{i_1}}$ such that
\begin{enumerate}
\item[(i)] $<x_1>\subset (\overline{J}_{i_1})_{y_{i_1}}.$
 
\item[(ii)] $(\call_{i_1}){y_{i_1}}=x_{j_1}^{c_1}\cdot\ldots\cdot x_{j_r}^{c_r}$,
with $1<j_1<\ldots<j_r<n$.
\end{enumerate}
 
If in addition  $R(1)(\Max t_{i_1}^d) = \emptyset$
locally at $y_{i_1}$ then  the basic
object $(Z_{i_1},(C(J^{\prime \prime}_{i_1},b^{\prime \prime}),
\overline{E}_{i_1}^{\prime \prime})$
can be defined by setting $Z_{i_1}=V(<x_1>)$, 
$C(J^{\prime \prime}_{i_1})$ as the
trace of $\overline{J}_{i_1}$ at ${\mathcal O}_{W_{i_1},y_{i_1}}/<x_1>$ 
and $\overline{E}_{i_1}^{\prime \prime}$ as the restriction of ${E}_{i_1}^{\prime \prime}$ to 
$Z_{i_1}$ 
(see 
\cite[9.3]{EncVil97:Tirol}).
 
\end{remark}

\section{Examples}
\label{examples}

The purpose of this section is to illustrate with some examples 
how the functions $f_i^d$ and $t_i^d$ are defined.

\begin{example}\label{ejemplo1} Let us consider again $W=\ma_{\rac}^3=\mbox{Spec }{\mathbb
{Q}}[x_1,x_2,x_3]$,  and let $X$ be the irreducible curve
determined by the ideal $J=<x_1,x_2x_3+x_2^3+x_3^3>=\cali(X)$.
 
\
 
Consider the basic object
$(W,(J,1),E=\{\emptyset\})=(W_0,(J_0,1),E_0)$ and set $f_0=x_2x_3+x_2^3+x_3^3$. 
Note that $\Sing(J_0,1)=X$ and that 
the order of the ideal $J_0$ at every point is 1. Since $E_0=\{\emptyset\}$
we see that $t^3_0(\xi)=(1,0)$ for any $\xi\in \Sing(J,1)$, and hence
$\max t^3_0=(1,0)$ and $\Max t^3_0=\Sing(J,1)=V(J)$.
 
\
 
Now we attach to $\Max t^3_0$ a basic object
$(W_0,(J_0^{\prime\prime},b^{\prime\prime}),E_0^{\prime\prime})$
in the sense of step 1 of  (\ref{secinduction}). In this case by \ref{nuestroalgo}, 
$J_0^{\prime\prime}=J_0$, 
$b^{\prime\prime}=1$ and $E_0^{\prime\prime}=\{\emptyset\}$.
 
\
 
Note that $R(1)(\Max
t^3_0)=R(1)(\Sing(J,1))=\emptyset$ since $\Sing(J,1)$ is a curve in a
3 space, so we proceed by induction in
the sense of Step 2 of \ref{secinduction}. Following \ref{nuestroalgo} we choose
the hypersurface $Z\subset W$ determined by the ideal $\cali(Z)=<x_1>$, and define 
\begin{equation}
\label{aghnadido}
C(J^{\prime\prime})=<x_2x_3+x_2^3+x_3^3> \subset {\mathbb
{Q}}[x_2,x_3]\simeq {\mathbb {Q}}[x_1,x_2,x_3]/<x_1>,
\end{equation}
and
$b^{\prime\prime}!=1$. So $Z=\ma_{\rac}^2\subset \ma_{\rac}^3=W$.
 
\
 
Recall that a resolution of the basic object
$(W_0,(J_0^{\prime\prime},b^{\prime\prime}),E_0)$ is equivalent to
a resolution  of
$(Z,(C(J^{\prime\prime}),b^{\prime\prime}!),E_Z)$. 

\
 
Now we define the function
$t^2_0:\Sing(C(J^{\prime\prime}),1)\longrightarrow
\rac\times\ent$ by 
        \begin{equation}
        \label{exp0}
t^2_0(\eta)=\left\{
        \begin{array}
        {lll}
         (2,0) & \mbox{ if } & \eta=(\bf{0},\bf{0})\in \ma_{\rac}^2\\
         (1,0) & \mbox{ if } & \eta\in X
\setminus\{(\bf{0},\bf{0})\}
        \end{array}
        \right.
        \end{equation}

Hence $\Max t^2_0=(\bf{0},\bf{0})$. Since $\Max t^2_0$ is a closed point and
the resolution function described in \ref{secinduction} is defined in terms of the functions 
 $t^j_i$, it
follows
that we blow-up at this closed point,
        \begin{equation}
        \label{exp1}
        (Z_0,(<f_0>,1),(E_Z)_0)\longleftarrow
        (Z_1,(<f_1>,1),(E_Z)_1=\{\overline{H}_1\}),
        \end{equation}
where $\overline{H}_1$ is the exceptional divisor and $f_1$ is the strict
transform of $f_0$. Note that
\begin{equation}
\label{exxp1}
<f_0>=\cali(\overline{H}_1)\cdot(<{f}_1>).
\end{equation}
 
Sequence (\ref{exp1}) induces a sequence
        \begin{equation}
        \label{exp11}
        (W_0,(J_0,1),E=\{\emptyset\})\longleftarrow
        (W_1,(J_1,1),E_1=\{H_1\}),
        \end{equation}
and in this case $J_1=\overline{J}_1$ (in the sense of
 (\ref{aaa1})).
 
\
 
Let $X_1\subset W_1$ be the strict transform of $X\subset W$,
and let $\{p,q\}=X_1\cap H_1$. Now we will describe locally 
$J_0\calo_{W_1,p}$: There is a regular system of parameters
$\{\frac{x_1}{x_2},x_2\}$ at $\calo_{W_1,p}$, such that
$J_0\calo_{W_1,p}=<\frac{x_1}{x_2},x_2(f_1)>$, and hence
\begin{equation}
\label{ninguna}J_0\calo_{W_1,p}\subset <\frac{x_1}{x_2},x_2>.
\end{equation}
Let $L_1=V(<\frac{x_1}{x_2},x_2>)$. Then $L_1\subset H_1$
and $\cali(L_1)$ will be an embedded component of the total
transform of $J_0$ at $\calo_{W_1}$. In particular this embedded
component will arise after any sequence of quadratic
transformations, and hence Hironaka's desingularization does not
fulfill property (iii) of Theorem \ref{MainTh}.
 
\
 
Following the resolution algorithm proposed in
\cite{EncVil99} (see also the addendum of \cite{EncVil97:Tirol}), the
desingularization  comes to an end  after blowing up at
the closed points $\{p,q\}$, and therefore the same argument as before
shows that property (iii) of Theorem \ref{MainTh} does not hold for
that algorithm.
\end{example}
 
\begin{example}
\label{ejemplo2} Now consider $W=\ma_{\rac}^4=\mbox{Spec }{\mathbb
{Q}}[x_0,x_1,x_2,x_3]$ and let $X$ be the irreducible closed
subscheme determined by the ideal
$J=<x_0,x_1,x_2x_3+x_2^3+x_3^3>=\cali(X)$. Let
$(W,(J,1),E=\{\emptyset\})=(W_0,(J_0,1),E_0)$ be the
corresponding basic object. Set $f_0=x_2x_3+x_2^3+x_3^3$. Once more
the function $t^4_0$ is constant on  $\Sing(J_0,1)$ with value
$(1,0)$, and hence $\max t^4_0=(1,0)$ and $\Max t^4_0=\Sing(J_0,1)=X$.
 
\
 
As in Step 1 of \ref{secinduction}, we associate to $\Max t^4_0=\Sing(J_0,1)$
the basic object $(W,(J^{\prime\prime}, b^{\prime\prime}),E),$  
where now 
$J^{\prime\prime}=J$, $b^{\prime\prime}=1$ and $E=\{\emptyset\}$ (see \ref{nuestroalgo}). 
So a resolution
of this basic object is a resolution of  $(W_0,(J_0,b),E_0)$.
 
\
 
Note that
$R(1)(\Sing(J^{\prime\prime},b^{\prime\prime}))=\emptyset$, so we
proceed by induction as in  Step 2 of  \ref{secinduction}. Let $Z$
be the regular hypersurface determined by the ideal
$\cali(Z)=<x_0>$, let 
\begin{equation}
\label{agnadido1} C(J^{\prime\prime})=<x_1,f>\subset
{\mathbb{Q}}[x_1,x_2,x_3]\simeq{\mathbb {Q}}[x_0,x_1,x_2,x_3]/<x_0>
\end{equation}
and $b^{\prime\prime}!=1$. Then a  resolution of the basic object
$(Z,(C(J^{\prime\prime}),b^{\prime\prime}!),(E_Z))$ is equivalent to a resolution 
 of the basic object 
$(W,(J^{\prime\prime},b^{\prime\prime}),E)$.

\
 
Note that the function  $t^3_0$ is constant along $\Sing(C(J^{\prime\prime}),1)$
and its value is $(1,0)$, 
hence 
$$\Max
t^3_0=\Sing(C(J^{\prime\prime}),1).$$

The next step of the algorithm of resolution is to associate a
basic object,
$$(Z^{\prime\prime},(C(J^{\prime\prime})^{\prime\prime},e),E_Z),$$
 to
$(Z,(C(J^{\prime\prime}),1),E_Z)$ in the sense of Step 1 of \ref{secinduction}.
Note that at this point we are in the same situation as in Example
\ref{ejemplo1}, thus a similar argument will tell us that we will
define a function $t^2_0$ in a two dimensional ambient space with
constant value equal to $(1,0)$.
 
\
 
Summarizing we have constructed the first three coordinates of
the function of the algorithm of
resolution at the first stage:
        \begin{equation}
        \label{ej2}
        \begin{array}
        {cccc}
        f^4_0: & X_0 & \longrightarrow & (\rac\times\ent)\times
        (\rac\times\ent)\times (\rac\times\ent)
         \end{array}
         \end{equation}
where,
$$\xi   \longrightarrow  f^4_0(\xi)=\left\{
\begin{array}{lll}
[(1,0),(1,0),(2,0)] & \mbox{ if } & \xi=(\bf{0},\bf{0},\bf{0},\bf{0})\in
\ma_{\rac}^4
\\ \left[(1,0),(1,0),(1,0)\right] & \mbox{ if } & \xi\in X\setminus
 (\bf{0},\bf{0},\bf{0},\bf{0})
\end{array}
\right.$$

\end{example}
 
\begin{example}
\label{ejemplo3} {\em Canonical choice of the centers.} Let
$W=\ma_{\rac}^3=\mbox{Spec }\rac[x_1,x_2,x_3]$,
\begin{equation}
\label{agnadido}J=<x_1,x_2>
\end{equation}
and
consider the associated basic object $(W,(J,1),E=\{\emptyset\})$.
The function $t^3_0$ is constant on $V(J)$ with value $(1,0)$ and
hence $\Max t^3_0=\Sing(J,1)$.
 
\
 
We associate the basic object
$(W,(J^{\prime\prime},b^{\prime\prime}),E=\{\emptyset\})$ to $\Max
t^3_0$ in the sense of Step 1 of \ref{secinduction}, where, once more,
 $J^{\prime\prime}=J$ and
$b^{\prime\prime}=1$.
 
\
 
As $R(1)(\Sing(J^{\prime\prime},b^{\prime\prime}))=\emptyset$, we
are under the assumptions  Step 2 of \ref{secinduction}, and hence we define
$Z=V(<x_1>)$,
$C(J^{\prime\prime})=<x_2>\subset \rac[x_2,x_3]\simeq
\rac[x_1,x_2,x_3]/<x_1>$, and $b^{\prime\prime}=1$.
 
\
 
Now we consider the basic object $(Z,(C(J^{\prime\prime}),1),
(E_Z)_0)$. Note that 
$$\Max t^2_0=R(1)(\Max
t^2_0)=\Sing(C(J^{\prime\prime}),1)=V(x_1,x_2)\subset
\ma_{\rac}^3,$$ 
and this is the canonical choice of center mentioned
in {\bf (1)} of step 2 in \ref{secinduction}.
\end{example}

\begin{remark}
\label{ejemplo31}
 Let us assume that $J\calo_W$ is an ideal
defining a closed regular and pure dimensional scheme $X$ of
codimension $e$ and consider the basic object $(W,(J,1),E=\emptyset)$. Then $\max
t^d=(1,0)$, $\Max t^d=X$ and
$(J^{\prime\prime}, b^{\prime\prime})=(J,1)$ as in Remark \ref{nuestroalgo}.
Locally at any point $\xi\in X$, the basic object $(Z,
(C(J^{\prime\prime}), l), E_{Z})$ is equal to
$(Z,
(J|_{Z}, 1), \emptyset)$. Now we are dealing with a
$d-1-$dimensional basic object where $V(J|_{Z})$ determines a
 regular scheme of codimension equal to the codimension of $X$ in $W$ 
(see Remark \ref{nuestroalgo}). So 
$(Z, (C(J^{\prime\prime}), l),\emptyset)$ is under the
same conditions as $(W,(J,1),E=\emptyset)$.
Therefore the first $e$ coordinates
of the function $f^d$ that defines the algorithm are:
\begin{equation}
\label{ultimoex}
\begin{array}{cccc}
f^{d,e}: & X & \longrightarrow &
\overbrace{(\rac\times\ent)\times\ldots\times
(\rac\times\ent)}^{e-\mbox{times}}\\
 & \xi & \longrightarrow & (1,0)\times\ldots\times (1,0)
\end{array}
\end{equation}
\end{remark}

\section{Proof of  Lemma \ref{lema}}
\label{last}
\begin{definition}
\label{rsing} 
Let  $X\subset W$ be under the
assumptions of Theorem  \ref{MainTh}, let 
$W_0=W$, and 
$J_0=J=\cali(X)$. Consider the basic object 
 $(W_0,(J_0,1),E_0)$  (note that in this case 
$\Sing(J_0,1)=X$). With this notation  define 
$$
\RSing(X)=\left\{y\in \Sing(J_0,b): 
\begin{array}{c}
y\in H_i \mbox{ for some } E_i\in E_0,
\mbox{ or }
J_0 \mbox{ is not regular at y}
\end{array}
\right\}
$$
The set $\RSing(X)$ will be called the {\em relative singular
locus} since the notion is relative to $E$. Note that  $\RSing(X)\subset X$ 
is closed. If $E_0=\emptyset$ then $\RSing(X)=\Sing(X)$. Although in the formulation of Theorem 
\ref{MainTh} we assume $E_0=\emptyset$, in our proofs, based in inductive arguments and step by step procedures, we will 
have to consider the case when $E_0\neq \emptyset$. 
\end{definition}
\begin{definition}
\label{property} [{\em The Relative Property.}]\hspace{0.2cm}{\rm 
With the same notation and assumptions as in Definition \ref{rsing}, consider
a sequence of monoidal transformations, 
\begin{equation}
\label{rsingpro}
(W_0,(J_0,1),E_0)\longleftarrow \ldots\longleftarrow (W_N,(J_N,1),E_N)
\end{equation}
at permissible centers $C_i\subset W_i$, and let $F \subset W_0$ be the union
of all images of $C_i$ (which will be a closed subset of $W_0$ since all
morphisms
$W_{i}\to W_0 $ are proper).  We
say that the sequence (\ref{rsingpro}) has the {\em relative property},
if\begin{equation}
        \label{relpro}
        F \subset \RSing(X).
        \end{equation}
        }
\end{definition}

\begin{definition}
\label{componentsr}
Given a sheaf of ideals $J\subset {\mathcal O}_{W}$, we will denote by 
$R(a)(J)$ the union of components of $V(J)$ of codimension $a$. Note that 
$R(a)(J)$ may be empty.  
\end{definition}

\begin{center}
{\bf Proof of Lemma \ref{lema}}
\end{center}

Let $(W_0,(J_0,1),E_0)$ be as in Theorem \ref{MainTh}. Note that every sequence of transformations of pairs 
$$(W_0,E_0)\longleftarrow \ldots\longleftarrow (W_N,E_N)$$
induces a  factorization $J_0{\mathcal O}_{W_N}=\call_N 
{\overline J}_N$. We will prove Lemma \ref{lema}
by induction on the relative codimension of the ideal $\overline{J}_N$. 
The two steps of the inductive argument are 
presented in \ref{inductionarg} and \ref{grande}.  

\begin{parrafo}
\label{inductionarg} {\rm {\bf Case 0.} {\em Let $X\subset W$ be 
under the assumptions of Theorem  \ref{MainTh}, and let 
$(W_{0},(J_{0},1),E_{0})=(W,(\cali(X),1),\emptyset)$. Then  
there is a finite sequence of
transformations of basic objects \begin{equation}
        \label{444}
        (W_{0},(J_{0},1),E_{0})
\longleftarrow\cdots\longleftarrow (W_{r_1},(J_{r_1},1),E_{r_1})
        \end{equation}
so that $J_{r_1}{\mathcal O}_{W_{r_1}}=\overline{J}_{r_1}$, and
$\overline{J}_{r_1}$  is of $(W_{r_1},E_{r_1})-$codimension
$\geq $ 1.}

\

The sequence (\ref{444}) will be defined   as a concatenation of two chains
of transformations: sequences (\ref{111}) and (\ref{tec2}), which will be
constructed, respectively, in two steps,  {\bf A} and {\bf B}. The condition
of relative codimension $\geq 1$ will be achieved really in Step {\bf A},
but Step {\bf B} will play a role in our inductive argument.

\
 
{\bf Step A.}
 
\
 
There is  a
sequence of transformations at permissible centers
$C_i\subset \Sing(J_1,1)$,  
\begin{equation}
        \label{111}
        (W_{0},(J_{0},1),E_{0})
\longleftarrow\cdots\longleftarrow (W_{k_1},(J_{k_1},1),E_{k_1})
        \end{equation}
so that
$$\max t^d_0 \geq \max t^d_1 \geq \max t^d_2 \ldots \geq
\max t^d_{k_1-1} > \max t^d_{k_1}$$ and $\max
t^d_{k_1}=(1,0)$ but $\max t_i^d > (1,0)$ for any index
$i<k_1$. If $\max t^d_0=(1,0)$ we would take $k_1=0$. 
Then,
\begin{equation}
        \label{555}
J_{k_1}=I(H_{1})^{a_{1}}\cdots I(H_{k_1})^{a_{k_1}}\overline{J}_{k_1}
        \end{equation}
where $a_1,\ldots,a_{k_1}\in \nat$, and by Remark \ref{svive}, 
$\overline{J}_{k_1}$  is of  $(W_{k_1},E_{k_1})-$codimension
$\geq 1$.

\

Let $U_0=W_0-\mbox{RSing}(X)$. Note that $U_0$ is dense in $W_0$, and that
$X\cap U_0$ is a smooth subscheme of $U_0$. 
Define 
$$
U_{k_1}\subset W_{k_1}
$$ 
as the pull-back of $U_0$ in  $ W_{k_1}$,
let $F_{k_1}$ be the union of the images of $C_i$ in  $W_0$, let 
$$V_0=W_0\setminus F_{k_1}$$  
and set $$V_{k_1}=W_{k_1}\setminus \cup_{H_i\in
E_{k_1}}H_i.$$ 

\
 
{\bf Step B.}
 
\
 
In expression (\ref{555}) consider the smallest index
$j$ so that $a_j \geq 1$. Since
$a_j \geq 1$, then $H_j \subset \Sing (J_{k_1},1)$, and it is also
clear that $H_j$ has normal crossings with $E_{k_1}$. Blow-up at $H_j$:
$$(W_{k_1},E_{k_1}) \longleftarrow (W_{k_1+1},E_{k_1+1}).$$
Then $W_{k_1}=W_{k_1+1}$, $E_{k_1+1}=E_{k_1}$ and       \begin{equation}
        \label{tec1}
        J_{k_1+1}=I(H_{1})^{a^*_{1}} \cdots
I(H_{i})^{a^*_{r}}\overline{J}_{k_1}
        \end{equation}
where $ a^*_{i}=a_{i}$ for each $i\in\{1,\ldots,k_1\}$, $i\neq j$,  and $
a^*_{j}=a_{j}-1$ (see (\ref{aaa1}) for the way that $J_{k_1+1}$ is defined). 
 
\
 
In this way we define, for some index $r_1 \geq k_1$, a sequence of
monoidal transformations :
        \begin{equation}
        \label{tec2}
        (W_{k_1},(J_{k_1},1),E_{k_1})
\longleftarrow\cdots\longleftarrow (W_{r_1},(J_{r_1},1),E_{r_1}),
        \end{equation}
so that
        \begin{equation}
        \label{tec3}
        (W_{k_1},E_{k_1})=(W_{r_1},E_{r_1})
        \end{equation}
but now expression (\ref{555}) reduces to
        \begin{equation}
        \label{tec4}
        J_{r_1}=\overline{J}_{r_1}
        \end{equation}
 
Note that $W_{r_1}\to W_{k_1}$ is the identity map. Set $ V_{r_1}$ and
$ U_{r_1}$ as pull-backs of $ V_{k_1}$ and $ U_{k_1}$ at $W_{r_1}$.
 
\
 
\begin{remark}
\label{remark0}
$\mbox{ }$ 
\begin{enumerate}
\item[i.] Note that $\max t^d_{r_1}=(1,0)$, thus
$\Max t_{r_1}^d=V(\overline{J}_{r_1})$ (so $t^d_{r_1}(y)=(1,0)$ for every 
$y \in V(\overline{J}_{r_1}))$, and therefore   by Remark \ref{svive} 
$\overline{J}_{r_1}$ is of $(W_{r_1},E_{r_1}^{\prime})-$codimension
$\geq 1$.
\item[ii.] The restriction of $\overline{J}_{r_1}(\subset {{\mathcal
O}}_{W_{r_1}})$ to $ V_{r_1}$ (respectively to $ U_{r_1}$) is naturally
identified with $J_0 $ restricted to $V_0$ (respectively with $J_0 $
restricted
to $U_0$).
\item[iii.] Since Definition \ref{omni} holds for $\overline{J}_{r_1}$ and
$(W_{r_1},E_{r_1})$ with $a=1$,  at any closed point
$y\in V(\overline{J}_{r_1})$, there is a
regular system of parameters $\{x_1,x_2,\ldots,x_d\}\subset
\calo_{W_{r_1},y}$ such that
$$<x_1>\subset (\overline{J}_{r_1})_y$$ 
and  if $y\in H_i\in E_{r_1}$
then $\cali(H_i)=<x_{i_j}>$ with $i_j>1$. This leads to:
\item[iv.] If $X$ is a subscheme of codimension $>1$ in $W_0$, then
$R(1)( V(\overline{J}_{r_1}))$  must be empty. In
fact, note that  by (iii), $R(1)(V(\overline{J}_{r_1}))\cap V_{r_1}$ must be dense in
$R(1)( V(\overline{J}_{r_1}))$; but $V_{r_1}$ is isomorphic to an open
subset of $W_0$. The claim follows now from (ii), and the
assumption on the codimension of $X$.
\end{enumerate}
\end{remark}
 
\noindent {\bf Claim:} {\em The concatenation of sequences
(\ref{111}) and (\ref{tec2}):
        \begin{equation}
        \label{tec5}
        (W_{0},(J_{0},1),E_{0}) \longleftarrow\cdots\longleftarrow
(W_{k_1},(J_{k_1},1),E_{k_1})\longleftarrow\cdots\longleftarrow
(W_{r_1},(J_{r_1},1),E_{r_1})
        \end{equation}
has the  relative property introduced in \ref{property}}.
 
\
 
\noindent{\em Proof of the Claim:} It suffices to check, that for
each index $i=1,\ldots,r_1$, the proper morphism $W_i \to W_0$ induces an
isomorphism over $U_0=W_0-\RSing(X)$.
 
\
 
For the construction in Step A, we may assume by induction that
$W_{i-1} \to W_0$ induces an isomorphism over $U_0=W_0-\RSing(X)$,
and note that $C_{i-1} \subset \Max t^d_{i-1}$. Since $\max
t^d_{i-1}>(1,0)$ it suffices to recall that $ t^d_{0}(y)=(1,0)$ for
any point $y\in \Sing(J_0,1)\cap U_0$.
 
\
 
The same statement of the claim  is clear from the construction in Step B,
since only exceptional hypersurfaces  are chosen as centers. This proves the
claim.}
 
\end{parrafo}
 
\begin{parrafo}
\label{grande}{\rm {\bf Case $e\geq 1$.} {\em Assume that $X$ is of
codimension $> e$ in $W_0$, and that  there is a  finite sequence of
transformations
of basic objects\begin{equation}
        \label{4444}
        (W_{0},(J_{0},1),E_{0})
{\longleftarrow}\cdots {\longleftarrow} (W_{k_e},(J_{k_e},1),E_{k_e})
        \end{equation}
at permissible centers $C_i$ so that:
 
\begin{enumerate}
 
\item[{\bf (1)}] If $F_{r_e}$ is the union of the images of $C_i$ in
$W_0$, if  $V_0^e=W_0\setminus F_{r_e}$,   and if $V_{r_e}\subset W_{r_e}$ is the
pull-back  of $V_0^e$ in $W_{r_e}$,   then  
$$ W_{0}\longleftarrow W_{r_e}$$
defines an isomorphism
$V^e_0\cong V_{r_e}$ and $F^e \subset \RSing(X)$ (i.e the relative
property in Definition \ref{property} holds); and hence $U_0\subset V^e_0$.
In particular, if $U_{r_e}\subset
W_{r_e}$ denotes the pull-back of $U_0$, then $U_{r_e}\subset
V_{r_e}$ and $U_{r_e}\cong U_0$.
\item[{\bf (2)}] $\overline{J}_{r_e}$ has relative codimension $\geq e$.
\item[{\bf (3)}] $R(e)(V(\overline{J}_{r_e}))= \emptyset$.
\end{enumerate}
Then under this assumptions there is an enlargement of the sequence (\ref{4444}), 
\begin{equation}
        \label{e1}
        (W_{0},(J_{0},1),E_{0}) \longleftarrow\cdots
\longleftarrow(W_{r_e},(J_{r_e},1),E_{r_e})\longleftarrow\cdots
\longleftarrow(W_{r_{e+1}},(J_{r_{e+1}},1),E_{r_{e+1}})
\end{equation}
so that $J_{r_{e+1}}=\overline{J}_{r_{e+1}}$   and
$\overline{J}_{r_{e+1}}$ has relative codimension   $\geq e+1$.}
 
\
 
We will accomplish this part of the proof in two steps, {\bf A} and {\bf B}.

\
 
{\bf Step A.}
 
\

Under the assumptions of \ref{grande},  we may 
assume that locally at $y \in V(J_{r_e})$:
\begin{enumerate}
\item[(i)]  All functions
$ t_{r_e}^{d}$,$ t_{r_e}^{d-1},\ldots, t_{r_e}^{d-e-1}$ are defined,
\item[(ii)] $\max t_{r_e}^d=(1,0),\max t_{r_e}^{d-1}=(1,0),\ldots,\max t_{r_e}^{d-e-1}=(1,0)$
\item[(iii)] By {\bf (2)} there is a regular system of parameters 
$\{x_1,\ldots,x_e\}\subset {\mathcal O}_{W_{r_e},y}$, so that:
\begin{enumerate}
\item[(a)] $<x_1,x_2, \ldots, x_e>\subset(\overline{J}_{r_e})_{y}\subset {\mathcal O}_{W_{r_e},y}$. 
\item[(b)] For each $H_i\in E_{r_e}$ with $y\in H_i$ there exists $j_i>e$ so that 
$$\cali(H_i)_y=<x_{j_i}>.$$
\end{enumerate}
\item[(iv)] By {\bf (3)}, $<x_1,x_2, \ldots , x_e> \neq (\overline{J}_{r_e})_{y}$.
\item[(v)] Since by {\bf (3)} $R(e)(V(\overline{J}_{r_e}))= \emptyset$, 
by Remark \ref{nuestroalgo}, 
\begin{equation}
\label{ccc}
\begin{array}{c}
(V(x_1,\ldots,x_e>),(\overline{J}_{r_e}|_{V(<x_1,\ldots,x_e>)},1),
(E_{r_e})|_{V(<x_1,\ldots,x_e>)})=\\
 \\
(\overline{W}_{r_e},({{\mathcal A}}_{r_e},1), \overline{E}_{r_e})
\end{array}
\end{equation}
is the $d-e$ dimensional basic object attached to $\Max t_{k_e}^{d-e}$ 
in  a suitable
neighborhood of $y\in {W}_{r_e}$. Note that ${\mathcal A}_{r_e}\neq 0$. 
This can be done for an open
covering, and these are the locally defined basic objects attached
to the value $$\underbrace{[(1,0),\ldots ,(1,0)]}_{e-\mbox{times}}$$
which are the first $e-$coordinates of the function$$f^d_{r_e}:
\Sing(J_{r_e},1)\longrightarrow \overbrace{T\times\ldots
\times T}^{d-\mbox{times}},$$
 (see \ref{secinduction}).
\end{enumerate}

\

By assumption we know that $J_{r_e}=\overline{J}_{r_e}$ is an ideal of
order 1 at every  point $y\in \Sing (J_{r_e},1)$. In  a  suitable open
neighborhood of $y$:
        \begin{equation}
        \label{ea1}
\Sing(J_{r_e},1)\subset
\overline{W}_{r_e}=V(<x_1,x_2, \ldots , x_e>)
        \end{equation}
and $\Sing(J_{r_e},1)= \Sing({{\mathcal A}}_{r_e},1)$ (at least locally).
 
\
 
Note that $\mbox{dim}(\overline{W}_{r_e})=(d-e)$, so $ t^{d-e}_{r_e}$ is
defined as a function on $\Sing({{\mathcal A}}_{r_e},1)$. 
By the previous identifications, we can view, at least locally, $
t^{d-e}_{r_e}$ as a function on $\Sing(J_{r_e},1)$.
By Remark \ref{versa}  these locally defined functions
$t^{d-e}_{r_e}$ patch so as to define a function on all
$V(J_{r_e})=\Sing(J_{r_e},1)$. Furthermore, there is a sequence 
of transformations of basic objects 
        \begin{equation}
        \label{ea2}
        (W_{r_e},(J_{r_e},1),E_{r_e})\longleftarrow\cdots\longleftarrow
(W_{k_{e+1}},(J_{k_{e+1}},1),E_{k_{e+1}})
        \end{equation}
so that
        \begin{equation}
        \label{ea3}
J_{j}=\overline{J}_{j}
        \end{equation}
for each index $j= r_e, r_e+1,\ldots, k_{e+1}$, which locally
induces  the sequence:
        \begin{equation}
        \label{ea4}
        (\overline{W}_{r_e},({{\mathcal A}}_{r_e},1),
\overline{E}_{r_e})\longleftarrow\cdots\longleftarrow
(\overline{W}_{k_{e+1}},({{\mathcal A}}_{k_{e+1}},1), \overline{E}_{k_{e+1}})
        \end{equation}
such that
$$\max t^{d-e}_{r_e} \geq \max t^{d-e}_{r_e+1} \geq
\ldots \geq
\max t^{d-e}_{k_{e+1}-1} > \max t^{d-e}_{k_{e+1}}=(1,0)$$ And hence
        \begin{equation}
        \label{fff}
{{\mathcal A}}_{k_{e+1}}=\cali(H_{r_e+1})^{a_{1}}\cdots
\cali(H_{k_{e+1}})^{a_{k_{e+1}-r_e}}\overline{{{\mathcal A}}}_{k_{e+1}}
        \end{equation}
where now $\overline{{{\mathcal A}}}_{k_{e+1}}$ has relative codimension
$\geq1$ (same argument as in Step A of Case 0, see also part (i) of
Remark \ref{remarke}).
 
\
 
\noindent{\em On the local description of $J_{k_{e+1}}$}:
 
\
 
Fix $y \in V(J_{k_{e+1}})$. At ${{\mathcal O}}_{W_{k_{e+1}},y}$ there is a
regular system of parameters
$ \{x_1,x_2, \ldots , x_d\} $ and there are   ideals $M_{k_{e+1}}$ and ${{\mathcal D}}_{k_{e+1}}$,
so that
        \begin{equation}
        \label{fff1}
        ({J_{k_{e+1}}})_y=<x_1,x_2, \ldots ,
x_e>+M_{k_{e+1}}\cdot{{\mathcal D}}_{k_{e+1}}
        \end{equation}
 where:
\begin{enumerate}
\item[i.] ${{\mathcal D}}_{k_{e+1}}$ induces
${\overline{\mathcal A}}_{k_{e+1}}$
(modulo $<x_1,x_2, \ldots , x_e>$), and $M_{k_{e+1}}$ is a monomial that
induces $I(H_{e_r+1})^{a_{1}}\cdots I(H_{k_{e+1}})^{a_{k_{e+1}-r_e}}$
modulo $(<x_1,\ldots,x_e>)$
(see (\ref{fff})).
\item[ii.] $x_{e+1}\in {{\mathcal D}}_{k_{e+1}}$ and $M_{k_{e+1}}$ is a monomial in
coordinates $x_s$ involving only indices $s>e+1$ (i.e.
$x_{e+1}\cdot M_{k_{e+1}}\in (J_{k_{e+1}})_y$).
\end{enumerate}
\
 
{\bf Step B.}
 
\
 
Now  we get  rid of the monomial $M_{k_{e+1}}$:  In expression  (\ref{fff}), consider the smallest
index $j$ so that there is a chart for which $a_j \geq 1$, and blow
up the hypersurface $\overline{H}_j$. Since $a_j \geq 1$, then
$\overline{H}_j \subset \Sing (\overline{{{\mathcal A}}}_{k_{e+1}},1)$, and
it is clear that $\overline{H}_j$ has normal crossings with
$\overline{E}_{k_{e+1}}$. In other words, take the smallest
index $j$ so that
        \begin{equation}
        \label{global}
        \mbox{dim}(H_j \cap
V(J_{k_{e+1}}))=d-e-1
        \end{equation}
Consider the blowing up at $H_j$, $$(W_{k_{e+1}},E_{k_{e+1}})\longleftarrow
(W_{(k_{e+1})+1},E_{(k_{e+1})+1})$$ Note that
 $\overline{W}_{k_{e+1}}=\overline{W}_{(k_{e+1})+1}$,
$\overline{E}_{k_{e+1}}=\overline{E}_{(k_{e+1})+1}$ and  that
        \begin{equation}
        \label{global1}
        {{\mathcal A}}_{(k_{e+1})+1}=I(H_{r_e+1})^{a^*_{1}} \cdots
I(H_{k_{e+1}})^{a^*_{r}}\overline{{{\mathcal A}}}_{(k_{e+1})+1}
        \end{equation}
where $ a^*_{i}=a_{i}$ for each $i\in\{r_e+1,\ldots,k_{e+1}\}$,  $i\neq j$, and $
a^*_{j}=a_{j}-1$. Note also that
        \begin{equation}
        \label{global2}
J_{(k_{e+1})+1}=\overline{J}_{(k_{e+1})+1}
        \end{equation}
and that $\Sing(J_{(k_{e+1})+1},1)\subset \overline{W}_{(k_{e+1})+1}$
and $\Sing(J_{(k_{e+1})+1},1)= \Sing({{\mathcal
A}}_{(k_{e+1})+1},1)$ (this follows from the properties listed
in  Step 2 of \ref{secinduction}).
 
\
 
We repeat this argument over  and over so as to define a sequence
 of
transformations of basic objects
        \begin{equation}
        \label{gr1}
        (W_{k_{e+1}},(J_{k_{e+1}},1),E_{k_{e+1}})\longleftarrow\cdots\longleftarrow
(W_{r_{e+1}},(J_{r_{e+1}},1),E_{r_{e+1}})
        \end{equation}
which induces a sequence
\begin{equation}
\label{grr2}
(\overline{W}_{k_{e+1}},(\cala_{k_{e+1}},1),\overline{E}_{k_{e+1}})
\longleftarrow\cdots\longleftarrow
(\overline{W}_{r_{e+1}},(\cala_{r_{e+1}},1),\overline{E}_{r_{e+1}})
\end{equation}
such that \begin{equation}
        \label{global3}
        {{\mathcal
A}}_{r_{e+1}}=\overline{{{\mathcal A}}}_{r_{e+1}}.
        \end{equation} 
Note that at each index $j$ of  the sequence (\ref{gr1}) :
        \begin{equation}
        \label{gr2}
        J_{j}=\overline{J}_{j}
        \end{equation}
and the local description in the last argument shows  that the
relative codimension of $J_{r_e+1}$ is $\geq e+1$.
 
\
 
\noindent {\bf Claim:} {\em The concatenation of the sequences (\ref{ea2})
and (\ref{gr1}):
        \begin{equation}
        \label{con1}
        (W_{0},(J_{0},1),E_{0})
\longleftarrow\cdots\longleftarrow
(W_{r_e},(J_{r_e},1),E_{r_e})\cdots\longleftarrow
(W_{r_{e+1}},(J_{r_{e+1}},1),E_{r_{e+1}})
        \end{equation}
has the relative property introduced in Definition \ref{property}}.
 
\
 
\noindent{\em Proof of the Claim:} It suffices to check, that for
each index $i$, the induced proper morphism $W_i \to W_0$ induces
an isomorphism over $U_0=W_0-\RSing(X)$.
 
\
 
For the construction in Step A, we argue by induction on $i\geq
r_e$, so assume that 
$$W_{i-1} \to W_0$$ 
induces an isomorphism over
$U_0=W_0-\RSing(X)$, and note that $C_{i-1} \subset \Max
t^{d-e}_{i-1}$. Since $\max t^{d-e}_{i-1}>(1,0)$ it suffices to
recall that $ t^{d-e}_{r_e}(y)=(1,0)$ for any point $y\in
\Sing(J_{r_e},1)\cap U_{r_e}$.
 
\
 
The statement of the claim  is clear from the construction in Step B, since
centers are chosen as exceptional hypersurfaces. This proves the
claim.}
\end{parrafo}
 
\begin{remark}
\label{remarke}
 Let  $U_{r_{e+1}}\subset W_{r_{e+1}}$ and
$V_{r_{e+1}}\subset W_{r_{e+1}}$ denote, respectively, the pull-backs of
$U_{r_e}\subset W_{r_e}$ and $V_{r_e}\subset W_{r_e}$ to $W_{r_{e+1}}$.
In particular,  $U_{r_{e+1}}\subset W_{r_{e+1}}$ denotes the pull-back of
$U_0$;
there is an inclusion $U_{r_{e+1}}\subset V_{r_{e+1}}$, and
$U_{r_{e+1}}\cong U_0$.

\begin{enumerate}
\item[i.] Since  $\max t^{d-e}_{r_{e+1}}=(1,0)$ we have that  
  $\Max t^{d-e}_{r_{e+1}}=V(\overline{J}_{r_{e+1}})$ (so
$t^{d-e}_{r_{e+1}}(y)=(1,0)$ for every $y \in V(\overline{J}_{r_{e+1}})$).
Set $J^{\prime\prime}=\overline{J}_{r_{e+1}}$,
$b^{\prime\prime}=1$ and $E^{\prime\prime}=E_{r_{e+1}}$ (in the sense
of \ref{secinduction}). Then $J^{\prime\prime}=\overline{J}_{r_{e+1}}$ is of
$(W_{r_{e+1}},E_{r_{e+1}})-$relative codimension
$\geq e+1$.
\item[ii.] The restriction of $\overline{J}_{r_{e+1}}\subset {{\mathcal
O}}_{W_{r_{e+1}}}$ to $ V_{r_{e+1}}$ (respectively to $ U_{r_{e+1}}$) is
naturally identified with $J_0 $ restricted to $V_0$ (respectively with $J_0 $
restricted to $U_0$).
\item[iii.]   Since Definition \ref{omni} holds for $\overline{J}_{r_{e+1}}$ 
 at every closed point $y\in V(\overline{J}_{r_{e+1}})$, there is a
regular system of parameters
$$\{x_1,\ldots,x_e,x_{e+1},\ldots,x_d\}\subset
\calo_{W_{r_e}}$$ 
such that
$<x_1,\ldots,x_e,x_{e+1}>\subset (\overline{J}_{r_{e+1}})_y$
and if
$y\in H_i\in E_{r_{e+1}}$, then $\cali(H_i)=<x_{i_j}>$ with $i_j>e+1$.
This leads to:
\item [iv.] If the codimension of $X$ in $W$ is $>e+1$, then
$R(e+1)( V(\overline{J}_{r_{e+1}}))$  must be
empty. In fact, (iii) asserts
that $R(e+1)( V(\overline{J}_{r_{e+1}}))\cap V_{r_{e+1}}$ must be dense in
$R(e+1)( V(\overline{J}_{r_{e+1}}))$; but $V_{r_{e+1}}$ is isomorphic to
an open subset of $W_0$. The claim follows now from (ii), and the
assumption that $X$ is of pure codimension $>e+1$. \qed
\end{enumerate}
\end{remark}

\end{document}